\documentclass{article}

\usepackage{arxiv}

\usepackage[utf8]{inputenc} 
\usepackage[T1]{fontenc}    
\usepackage{hyperref}       
\usepackage{url}            
\usepackage{booktabs}       
\usepackage{amsfonts}       
\usepackage{nicefrac}       
\usepackage{microtype}      
\usepackage{lipsum}
\usepackage{amsmath}
\usepackage{inputenc}

\usepackage{subfig}
\usepackage{rotating}
\usepackage{tikz}
\usepackage{lscape}
\usetikzlibrary{positioning}
\usepackage{array,calc}
\newlength\lena \newlength\lenb \newlength\lenc \newlength\lend
\newcolumntype{P}[1]{>{\centering\arraybackslash}p{#1}} 

\title{Deep Learning Based Anticipatory Multi-Objective Eco-Routing Strategies for Connected \& Automated Vehicles}

\author{
	\noindent
  Lama Alfaseeh \\
  Laboratory of Innovations in Transportation (LiTrans)\\
  Ryerson University\\
  Toronto, Canada \\
  \texttt{lalfaseeh@ryerson.ca}
  \And   
  Bilal Farooq \\
  Laboratory of Innovations in Transportation (LiTrans)\\
  Ryerson University\\
  Toronto, Canada \\
  \texttt{bilal.farooq@ryerson.ca}
  \\
}

\begin{document}
\maketitle

\begin{abstract}
This study exploits the advancements in information and communication technology (ICT), connected and automated vehicles (CAVs), and sensing, to develop anticipatory multi-objective eco-routing strategies. For a robust application, several GHG costing approaches are examined. The predictive models for the link level traffic and emission states are developed using long short term memory deep network with exogenous predictors. It is found that anticipatory routing strategies outperformed the myopic strategies, regardless of the routing objective. Whether myopic or anticipatory, the multi-objective routing, with travel time and GHG minimization as objectives, outperformed the single objective routing strategies, causing a reduction in the average travel time (TT), average vehicle kilometre travelled (VKT), total GHG and total NOx by 17\%, 21\%, 18\%, and 20\%, respectively. Finally, the additional TT and VKT
experienced by the vehicles in the network contributed adversely to the amount of GHG and NOx produced in the network.
\end{abstract}




\section{Introduction}
\label{Introduction}
\textcolor{black}{Transportation system in the U.S. produced 28\% of the total greenhouse gas (GHG) emissions in 2018, which is the largest share from a single source \cite{epa2020sources}.} According to the United States Environmental Protection Agency in 2017, 81\% of the GHG was CO$_2$, which is a major contributor to climate change and global warming. The employment of information and communication technology (ICT) and connected and automated vehicles (CAVs) has been suggested as a potential solution to alleviate the undesirable social and environmental impacts of transportation systems \cite{zegeye2009model}. \textcolor{black}{Particularly, routing schemes for CAVs have an explicit impact on traffic \cite{Djavadian2018Distributed, alfaseeh2018impact, lee2008evaluation, yang2006modeling} and environmental \cite{tu2019quantifying} characteristics.}


\textcolor{black}{While myopic routing is reactive, the anticipatory routing is proactive \cite{bottom2000consistent}. Myopic routing uses the current network state, while the anticipatory routing exploits the predicted conditions, while routing vehicles to their destinations \cite{bottom2000consistent}. Anticipatory routing is considered a promising approach to improve traffic characteristics of a network, while avoiding congestion--\textcolor{black}{especially when we employ a high market penetration rate (MPR) of vehicles that are equipped with a routing system that is based on anticipatory information} \textcolor{black}{\cite{bottom2000consistent, ben2001network, mahmassani1994development}}.}

\textcolor{black}{Eco-routing is a special case of routing that specifically considers the environmental aspects \cite{luo2016real, tzeng1993multiobjective}.} \textcolor{black}{Several studies in the literature developed myopic eco-routing systems. \textcolor{black}{For instance, \cite{djavadian2020multi} and  \cite{alfaseeh2019multi} developed single as well as multi-objective eco-routing for CAVs and found encouraging results.} A review paper developed by \cite{alfaseeh2020multifactor} classified the eco-routing models and illustrated the strengths and weaknesses of each category. The common limitations are related to the \textcolor{black}{level of data point resolution}, scale of the case study, and number of objectives optimized at once.} \textcolor{black}{With reference to anticipatory routing in general, a scarcity of studies, especially recent ones, in the literature has been noticed.} For the few anticipatory routing studies considering only travel time as the routing objective \cite{pan2013proactive, ben2001network, bottom2000consistent}, the limitations \textcolor{black}{associated are related to the scale of the case study \cite{kaufman1991iterative}, level of temporal and spatial resolution \cite{bilali2019analyzing}, use of centralized solutions that suffer from scaling issues, and the use of reflective prediction models \cite{pan2013proactive, ben2001network, bottom2000consistent,kaufman1991iterative}.} The anticipatory routing studies found, did not employ \textcolor{black}{sophisticated} predictive models, which is a major limitation. The forecasted travel time data points  were provided by running the traffic simulation in advance \cite{bilali2019analyzing, kim2016dynamic,ben2001network} or based on regression models between speed and other traffic variables, such as density  \cite{pan2013proactive}. In most of the studies, travel time of time step $t+i$, obtained from the traffic simulation or historical data, was used in time step $t$ for the anticipatory routing application, where $i$ is the prediction interval considered \cite{ben2001network}.

\textcolor{black}{Unlike previous studies, this work develops anticipatory multi-objective eco-routing strategies that can be implemented in the routing systems for connected and automated vehicles. Similar to the myopic routing schemes proposed in \cite{djavadian2020multi}, a microscopic level of aggregation, an urban network as a case study, and a high level of spatial \textcolor{black}{(link level)} and temporal \textcolor{black}{(1 minute)} resolution are employed in this study. \textcolor{black}{Unlike the previous anticipatory approaches that were centralized solutions,} this study develops the anticipatory routing strategies based on a dynamic distributed routing framework i.e. End-to-End Routing for Connected and Automated Vehicles (E2ECAV) \cite{FarooqandDjavadian}. In this study three eco-routing strategies are applied to optimize travel time (TT), GHG, as well as the combination of TT  and GHG. The performance indicators considered for comparison between the different routing strategies are average TT, average vehicle kilometres travelled (VKT), total GHG, and total NOx emissions produced. \textcolor{black}{Following are the main contributions}:
	\begin{itemize}
		\item Examination of the most representative GHG costing approach at highly disaggregate spatial (link level) and temporal (1 minute prediction interval) resolution.
		\item Development of anticipatory multi-objective eco-routing strategies, while employing the predictive model for speed developed here and a GHG prediction model developed in \cite{Lama_Prediction2020}.
		\item Detailed comparison between myopic and anticipatory routing strategies.
\end{itemize}}

\textcolor{black}{A brief literature review related to the existing eco-routing studies, their strengths, and weaknesses is presented in Section \ref{Literature Review}. In addition, existing predictive models of travel time and GHG and anticipatory routing studies are illustrated in Chapter \ref{Literature Review}. The specifications of the case study deployed is in Section \ref{Case study}. Section \ref{Methodology} includes the details related to the traffic and emission models, GHG costing approaches, and GHG and speed predictive models. Section \ref{Discussion and results} incorporates the results and discussions related. Finally, Section \ref{Conclusion} summarizes the major findings and future outlook.}


\section{Background}
\label{Literature Review}
\textcolor{black}{Due to the} advancements in sensing and communication technologies, their utilization in transportation systems, and emergence of connected and automated vehicles, the availability of the real-time high-resolution data has become feasible at large-scale. \textcolor{black}{Such data are adopted to develop highly accurate prediction models that have the potential to be used for anticipatory routing with multiple objectives. In this section, the studies found in the literature related to the myopic eco-routing, predictive models of time and GHG, and anticipatory routing are briefly presented.}

\cite{alfaseeh2020multifactor} developed a comprehensive review of eco-routing studies. They reported that the previous myopic eco-routing studies have predominantly used macroscopic \textcolor{black}{level traffic and emission models, are based on small case studies}, used centralized routing mechanisms, and have optimized single routing objective at a time. To overcome the aforementioned limitations, \cite{alfaseeh2019multi} applied multi-objective eco-routing in a distributed routing \textcolor{black}{framework}. The authors used the per lane weighted average for the GHG cost on links. Although normalizing by the number of lanes resulted in the underestimation for the links with a large number of lanes, reductions in the travel time and emissions produced were noticed when the multi-objective routing was applied. \cite{djavadian2020multi} developed myopic multi-objective eco-routing strategies for CAVs, considering the cost of idling as a penalty at the downstream intersection of a link at every interval. \textcolor{black}{The authors used the marginal cost for GHG in the objective function.} They found that including the penalty cost contributed to  \textcolor{black}{reductions of 4\% and 3\% in the average travel time and total GHG produced, respectively, in the case of the multi-objective routing when compared to the single objective routing.}

The predictive models are an essential component of the anticipatory routing. \cite{vlahogianni2014short} conducted a comprehensive review related to the short-term traffic forecasting. It was noticed that studies mainly considered freeways as their case studies, statistical models, \textcolor{black}{and a temporal resolution of five minutes in most of the cases \cite{vlahogianni2014short}.
	Freeways are utilized due to the complexity associated with urban congested networks. The vehicular dynamics in an urban areas are subject to changes during a \textcolor{black}{short time period (seconds). It is the result of stop and go phenomena and shorter length of the links, when compared to the freeways.} A low level of temporal resolution is employed due to the scarcity of microscopic data points and the high computational power required.} \textcolor{black}{The statistical models, \textcolor{black}{for example autoregressive integrated moving average (ARIMA) model}, are easy to use, but have limitations when dealing with complex non-linear relationships between the variables concerned \cite{zhang2003time}.}

When it comes to the travel time prediction, there are two main streams. One that predicts travel time directly, while the other stream predicts speed and consequently travel time. For directly forecasting travel time, several predictive models have been employed. Linear modelling \cite{zhang2003short}, nonlinear autoregressive with external inputs (NARX) model \cite{mane2018link}, nonlinear autoregressive model (NAR) \cite{mane2018link}, clustering \cite{elhenawy2014dynamic}, neural networks (NNs) \cite{mane2018link}, and deep neural networks \cite{ran2019lstm,duan2016travel}. 

A large body of literature exists where travel time is implicitly predicted from speed,  including \cite{gu2019short,yao2017short,ma2015long,yildirimoglu2013experienced, ishak2004optimizing,innamaa2000short}. \textcolor{black}{The common features between most of the aforementioned predictive models are the low temporal resolution and small case study employed. Not to forget that the statistical models are dominating, despite their inability to capture the complicated relationship between the variables in concern.} Even when a large network was employed as in \cite{zhang2019link}, the speed was predicted at \textcolor{black}{a low temporal resolution level}. With regards to the predictive models for speed and travel time, it was found in the literature that \textcolor{black}{long short terms memory} (LSTM) outperformed other predictive approaches including the ARIMA model \cite{ma2015long}.

Related to the GHG predictive models, GHG emissions were predicted based on yearly data points of fuel \cite{zhao2011mapreduce}, gross domestic product, or other economical factors \cite{ameyaw2019investigating, ameyaw2018analyzing, pao2011modeling}. The predictive models varied from statistical \cite{rahman2017modeling, tudor2016predicting} to deep neural networks based \cite{ameyaw2019investigating}. \textcolor{black}{ To overcome the limitations of the previous predictive models, the low spatial (national) and temporal (year) resolution, \cite{Lama_Prediction2020} developed a predictive model based on LSTM. The GHG emission rate (ER) at a link level and one minute time resolution was predicted based on the most representative traffic indicators of the previous time intervals.} 

\textcolor{black}{With reference to the anticipatory routing, several frameworks were proposed \cite{ben2001network, bottom2000consistent, mahmassani1994development} to minimize travel time. \cite{bottom2000consistent} developed a framework for the anticipatory routing in a \textcolor{black}{network with only a single OD pair, 14 links, and eleven OD paths}. The author developed a simulator with three major components. \textcolor{black}{The three components} considered are all time-dependent and consist of network conditions, path splits, and guidance messages. \textcolor{black}{Three maps were used to illustrate the relationships between the aforementioned three variables. The network loading map, used the path splits to define the network conditions. The guidance map employed the network conditions to define the guidance messages. Finally, the routing map translated the guidance messages into path splits.} 
	\textcolor{black}{\textcolor{black}{In terms of the predicted traffic variables, the real-time traffic characteristics and other related data were forecasted at short- and medium-term for the anticipatory routing.}} The variable message signs (VMS) were employed to provide vehicles with the best route. When the best route is defined, the driver's compliance was not guaranteed. Hence, the author incorporated a logit model to define the drivers' path choice \cite{bottom2000consistent}. \cite{ben2001network} proposed DynaMIT, which \textcolor{black}{can be employed to generate real-time guidance provided to the drivers}. Off-line and real \textcolor{black}{time} information were adopted. The \textcolor{black}{off-line data as well as historical network conditions} that were used for the state estimation. While the real-time data were obtained from the control system. Two simulation tools were used, \textcolor{black}{a demand and supply simulator. The demand simulator} estimated and forecasted the origin-destination (OD) flow, departure time of drivers, model, and route choice. While the supply simulator directly simulated the interactions between the demand and supply (network). Anticipatory routing was applied \textcolor{black}{while travel time minimization was the \textcolor{black}{routing} objective.} The predicted travel time was a function of \textcolor{black}{experienced travel time from the previous iteration}. Speed on links was estimated based on the linear relationship with density on the link in concern. The VMS were used to inform drivers with the information of the best route. With regards to the findings, the authors found that anticipatory routing is promising as it contributed to reductions in the travel time of vehicles \cite{ben2001network}.} \cite{pan2013proactive} suggested proactive re-routing strategies to reduce travel time. They \textcolor{black}{predicted} congestion based on density/jam density ratio. The speed was \textcolor{black}{predicted} based on the Greenshield model, linear relationship between density and speed, which is a limitation \textcolor{black}{as realistically the relationship is not linear. When the density is low, the speed is underestimated \cite{papacostas1993transportation}}. The authors found that their rerouting performed as good as the dynamic traffic assignment \cite{pan2013proactive}. \textcolor{black}{Another example is the work by \cite{liang2014real}, who applied re-routing based} on congestion prediction. One of the predictive models developed was based on the spatiotemporal correlation. The authors assumed that the traffic \textcolor{black}{flow} was constant during each prediction time interval, which is unrealistic. \cite{liu2016dynamic} proposed \textcolor{black}{a dynamic congestion model based on crowdsourcing in order to apply} the anticipatory routing for \textcolor{black}{a set of cooperative vehicles} by predicting the \textcolor{black}{probability} distribution of traffic conditions. The time interval adopted for the routes update was 1 minute and the data was obtained from the GPS traces and social media. They found that their approach outperformed the \textcolor{black}{myopic} routing approach \cite{liu2016dynamic}.

\textcolor{black}{To summarize, the existing predictive models are associated with limitations related to the spatial and temporal resolution. The anticipatory routing studies did not adopt efficient and} \textcolor{black}{more accurate,} predictive models and \textcolor{black}{were used in the context of a centralized routing framework}. \textcolor{black}{The centralized routing frameworks require a large infrastructure investment, are highly sensitivity to system failures, and involve high degree of complexity in the case of a system upgrade} {\cite{FarooqandDjavadian}. Hence, this study will tackle the aforementioned limitations. To the best of our knowledge, our study is the first of its kind to apply the anticipatory multi-objective eco-routing while deploying deep learning based predictive models in a distributed routing framework.}

\section{Case study}
\label{Case study}
\begin{figure*}[ht]
	\centering
	\includegraphics[width=\textwidth]{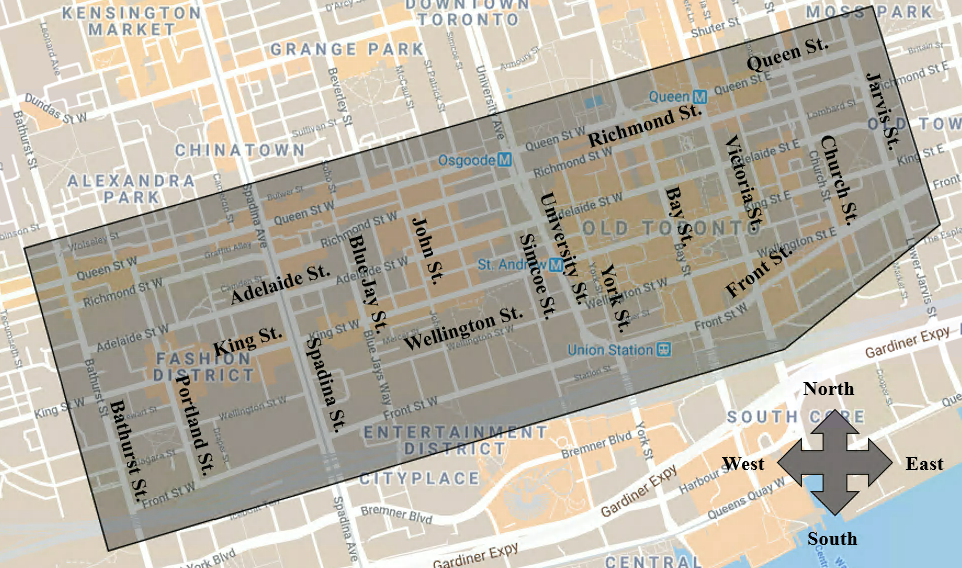}
	\caption{Case study, downtown Toronto}
	\label{Case_study}
\end{figure*}

\textcolor{black}{Downtown Toronto's road network is adopted as a case study because} it experiences high levels of recurrent congestion--especially during the morning peak period \textcolor{black}{\cite{force2012high}}. Downtown Toronto is the financial centre of Canada and has the highest job density among the major cities in the country. The network \textcolor{black}{is composed of} 223 links and 76 \textcolor{black}{intersections}. \textcolor{black}{Based on the 2019/2020 Toronto's vital report \cite{torontovitalsigns}, several factors contribute \textcolor{black}{to the excessive congestion levels in Toronto}. The population of Toronto has increased yearly by 1\% since 2011. Due to high cost of living, Toronto is the most expensive major city in the country. The ownership costs are growing four times faster than income, while renting costs are growing two times faster than income over the last decade \cite{torontovitalsigns}.} \textcolor{black}{The vehicular demand is provided by the Transportation Tomorrow Survey (TTS) for the period between 7:45am and 8:00am for the year 2014.} \textcolor{black}{Links in the case study are associated with different features with respect to the speed limit, number of lanes, and number of directions, i.e. a high level of heterogeneity is assured for a generic application.}
{Figure \ref{Case_study}} illustrates the area, including the major roads. The high level of heterogeneity assures a realistic generic application, especially for prediction. The speed limit on links is 2\%, 1\%, 30\%, 59\%, and 8\% of 10, 30, 40, 60, and 80 (km/h), respectively. With regards to the number of lanes, 1, 2, 3, and 4 of 7\%, 71\%, 15\%, and 7\%, respectively are used. 

\section{Methodology}
\label{Methodology}
\begin{figure*}[!ht]
	\centering
	\includegraphics[width=5in]{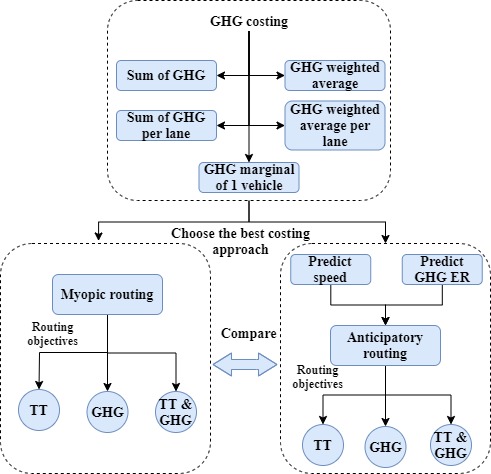}
	\caption{Methodology followed in this study}
	\label{Non_myopic_routing_methodology}
\end{figure*}

This section includes the specifications of the traffic and emission models, GHG costing approaches investigated, GHG and speed predictive models adopted, and the routing strategies taken into account. Figure \ref{Non_myopic_routing_methodology} demonstrates the general framework followed in this study.

\subsection{Traffic and emission models}
\label{Traffic and emission models}
\textcolor{black}{Microscopic traffic \cite{Djavadian2018Distributed} and emission simulators \cite{MOVES} are deployed in this study to obtain high resolution data points at every second. The Intelligent Driver Model (IDM) \cite{Treiber2000} is the car-following adopted} for the displacement estimation at every second, which is used to calculate the speed of vehicles \cite{Djavadian2018Distributed}. The second-by-second vehicular characteristics are captured and then used to estimate the space mean link indicators. \textcolor{black}{When all of the vehicles reach their destinations, the} simulation ends. \textcolor{black}{The indicators of links, speed, density, GHG, and flow, are updated at every} minute. 

The Motor Vehicle Emission Simulator (MOVES) is the emission model employed to estimate the second by second GHG (in $\text{CO}_{2\text{eq}}$) of every vehicle \cite{MOVES}. \textcolor{black}{The second-by-second emissions produced by vehicles are estimated based on} the vehicle operating mode, which \textcolor{black}{depends} on the vehicle specific power (VSP) \cite{MOVES}.

\subsection{GHG costing approaches}
\label{GHG costing approaches}

For a robust eco-routing strategy, five different GHG costing approaches have been assessed as illustrated in Table \ref{GHG_costing_strategies}. GHG$_{cost1(l,\Delta_j)}$, as in Equation \ref{GHG_cost1}, is for when GHG cost is the sum of GHG emissions of vehicles $N$ on the studied link $l$ during and interval $\Delta_j$. $\lambda_{(i,l,k)}$ is a binary variable, 1 if vehicle $i$ is on the studied link $l$ at time $k$ and 0 otherwise. \textcolor{black}{As the links in the network have different number of lanes,} GHG$_{cost2(l,\Delta_j)}$ normalizes the total GHG cost based on the number of lanes $Z_l$ \textcolor{black}{following} Equation \ref{GHG_cost2}. Considering a higher temporal resolution is illustrated in the GHG$_{cost3(l,\Delta_j)}$ costing approach. Weighted average of GHG produced on link $l$ is the outcome of Equation \ref{GHG_cost3}. At every second $k$ of any interval (1-60 \textcolor{black}{seconds}), a weight $k$, is multiplied by the GHG produced by the vehicles. \textcolor{black}{It means that the GHG emissions produced at $k$ of second 1 is multiplied by a weight of 1, while the GHG emission produced at second $k$ = 60 is multiplied by a weight of 60. The most recent seconds of an interval $\Delta_j$ are associated with a higher weight in the link cost.} The sum is then divided by the sum of weights. This costing approach takes the GHG cost at every second, which means it is associated with a high temporal resolution compared to GHG$_{cost1(l,\Delta_j)}$ and GHG$_{cost2(l,\Delta_j)}$. GHG$_{cost4(l,\Delta_j)}$ as in Equation \ref{GHG_cost4} follows the same logic of GHG$_{cost3(l,\Delta_j)}$, but divides by the number of lanes $Z_l$ of link $l$ to normalize. Finally, GHG$_{cost5(l,\Delta_j)}$ is the marginal cost of one vehicle traversing the studied link $l$. The marginal cost, as shown in Equation \ref{GHG_cost5}, depends on an estimated emission rate (ER) and the TT of interval $\Delta_j$ on link $l$. TT is obtained following Equation \ref{TT_cost}, where $D_l$ represents link $l$ length and $V(l,\Delta_j)$ is link $l$  speed of $\Delta_j$ time interval.

\begin{eqnarray}
\centering
\displaystyle
GHG_{cost1(l,\Delta_j)}=\sum_{k=1}^{\Delta_j}\sum_{i=1}^{N}{\lambda_{(i,l,k)} \times GHG_{(i,k)} }
\label{GHG_cost1}
\end{eqnarray}

\begin{eqnarray}
\centering
\displaystyle
GHG_{cost2(l,\Delta_j)}=\frac{\sum_{k=1}^{\Delta_j}(\sum_{i=1}^{N}{\lambda_{(i,l,k)} \times GHG_{(i,k)} })}{Z_l}
\label{GHG_cost2}
\end{eqnarray}

\begin{eqnarray}
\centering
\displaystyle
GHG_{cost3(l,\Delta_j)}=\frac{\sum_{k=1}^{\Delta_j} (k \times \sum_{i=1}^{N}{\lambda_{(i,l,k)} \times GHG_{(i,k)} })}{\sum_{k=1}^{\Delta_j}{k}}
\label{GHG_cost3}
\end{eqnarray}

\begin{eqnarray}
\centering
\displaystyle
GHG_{cost4(l,\Delta_j)}=\frac{\sum_{k=1}^{\Delta_j} (k \times \frac{ \sum_{i=1}^{N}{\lambda_{(i,l,k)} \times GHG_{(i,k)} }}{  Z_l})}{\sum_{k=1}^{\Delta_j}{k}}
\label{GHG_cost4}
\end{eqnarray}

\begin{equation}
GHG_{cost5(l,\Delta_j)}= GHG_{ER(l,\Delta_j)} \times TT_{(l,\Delta_j)}
\label{GHG_cost5}
\end{equation}

\begin{eqnarray}
TT_{(l, \Delta_j)}= \frac{D_l}{V(l, \Delta_j)}
\label{TT_cost}
\end{eqnarray}

\begin{table}[!h]
	\caption{GHG costing strategies investigated}
	\begin{center}
		\small
		\begin{tabular}[!t]{l  c c}
			\hline
			Costing strategy name & Approach & Equation used \\
			\hline
			\hline
			\centering
			GHG$_{cost1}$& Sum of GHG & \ref{GHG_cost1} \\
			\hline
			GHG$_{cost2}$& Sum of GHG per lane & \ref{GHG_cost2}\\
			\hline
			GHG$_{cost3}$ & Weighted average of GHG & \ref{GHG_cost3} \\
			\hline
			GHG$_{cost4}$& Weighted average of GHG per lane & \ref{GHG_cost4}\\
			\hline
			GHG$_{cost5}$& Marginal GHG cost of 1 vehicle & \ref{GHG_cost5}\\
			\hline
		\end{tabular}
	\end{center}
	\label{GHG_costing_strategies}
\end{table}

\subsection{GHG emission rate and speed predictive models}
\label{GHG emission rate and speed predictive models_methodology}

\textcolor{black}{Two separate LSTM networks have been trained to predict the variables, TT and GHG, required for the anticipatory multi-objective routing.} LSTM has been chosen as it outperforms the statistical models for time series data, such as ARIMA \cite{Lama_Prediction2020}. It also overcomes the shortcoming of standard neural networks, such as the vanishing gradient problem \cite{amarpuri2019prediction}. 

With reference to the LSTM architecture and hyper-parameters, it has been found that the selection of predictors, number of sequences, and the set of hyper-parameters \cite{reimers2017optimal, hutter2015beyond} have an explicit impact on the prediction performance. In addition, increasing the depth of the NNs may also introduce further enhancements \cite{pascanu2013construct, hermans2013training}. Not to forget that the efficiency of the hyper-parameters tuning process is profound \cite{snoek2012practical}. In this study, a comprehensive correlation analysis has been conducted for each of the predicted responses as discussed in Section \ref{Correlation analysis of the predicted variables}, GHG ER and speed. For each of the predicted variables, i.e. GHG ER and speed, first the \textcolor{black}{most representative predictors and number of previous time steps (sequences) used in the model are defined} based on the correlation analysis. Then hyper-parameters are tuned in two stages. The manual tuning mainly tries to narrow the search range of the hyper-parameters in concern for a more efficient systematic tuning process based on the Bayesian optimization \cite{wu2019hyperparameter}. Further details related to the predictive LSTM network can be found in \cite{Lama_Prediction2020}. To compare between the trained LSTM networks, four indicators are \textcolor{black}{utilized: 1) correlation coefficient between observed and predicted GHG ERs (in $\text{CO}_{2\text{eq}}$ g/sec), 2) fit to the ideal 45$^o$ line, 3) root mean square error (RMSE), and 4) R$^2$ }.

\subsubsection{Data collection}
\label{Data collection}
This step is essential for the development of the predictive models of GHG ER and speed. The quality of the data collected contributes to how reflective the predictive models are. The data points are extracted from an agent-based traffic model developed in \cite{Djavadian2018Distributed}. The demand is synthesized based on actual data from the Transportation Tomorrow Survey (TTS). \textcolor{black}{To train the LSTM networks, 80\% of the data is employed, while 20\% is used for testing.} The training and testing sets are 48,652 and 12,159 data points, respectively for the LSTM predictive models. The high level of heterogeneity of the traffic and environmental variables contributes to more generic predictive models. A wide \textcolor{black}{variety} of traffic and environmental conditions are captured and used for training the predictive models. To produce different traffic conditions, different demand levels and different departure time distributions are adopted. 



\begin{figure}[!h]
	\centering
	\subfloat[]{\includegraphics[width=2.5in]{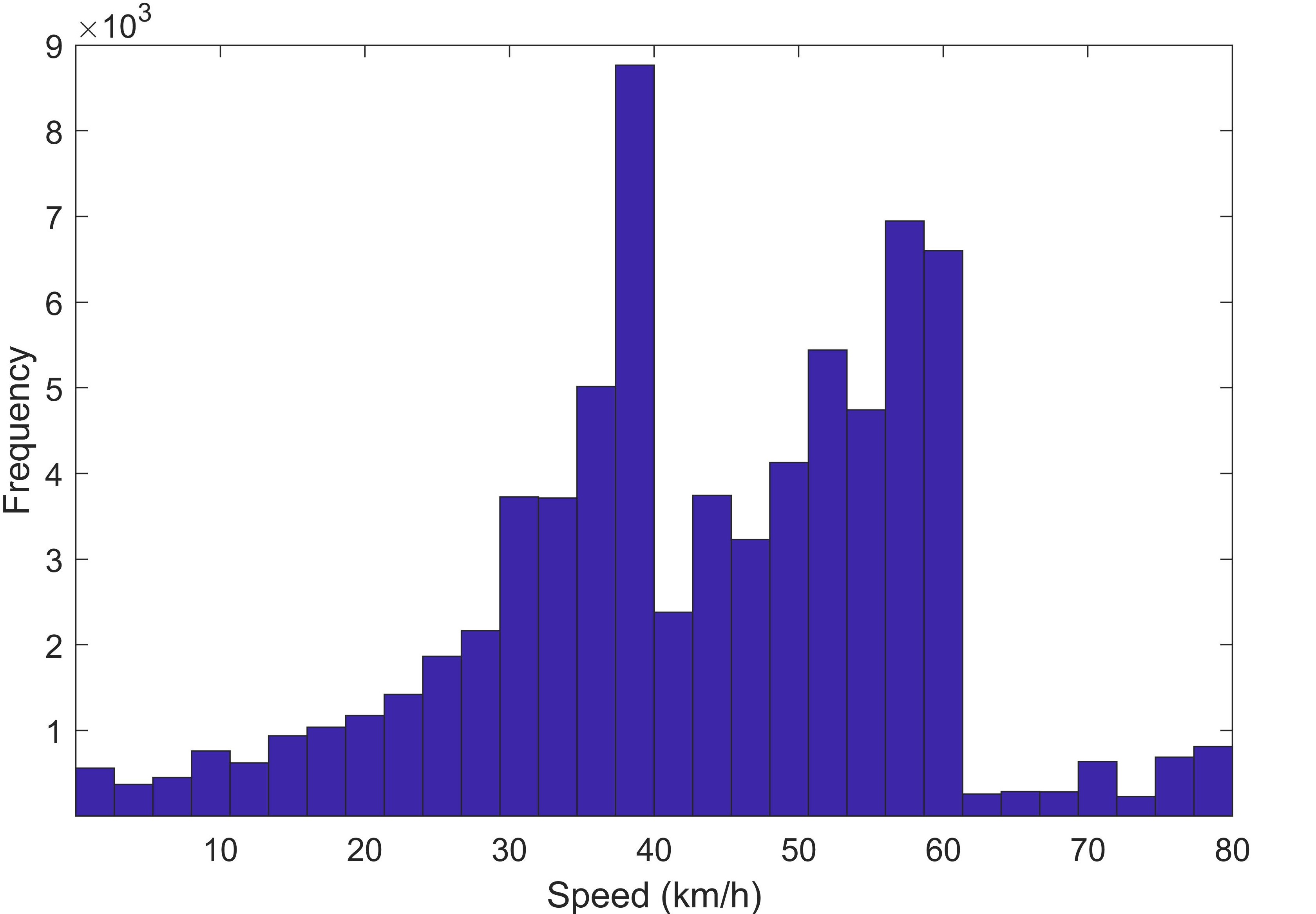}
		\label{histogram speed}}
	\hfil
	\subfloat[]{\includegraphics[width=2.5in]{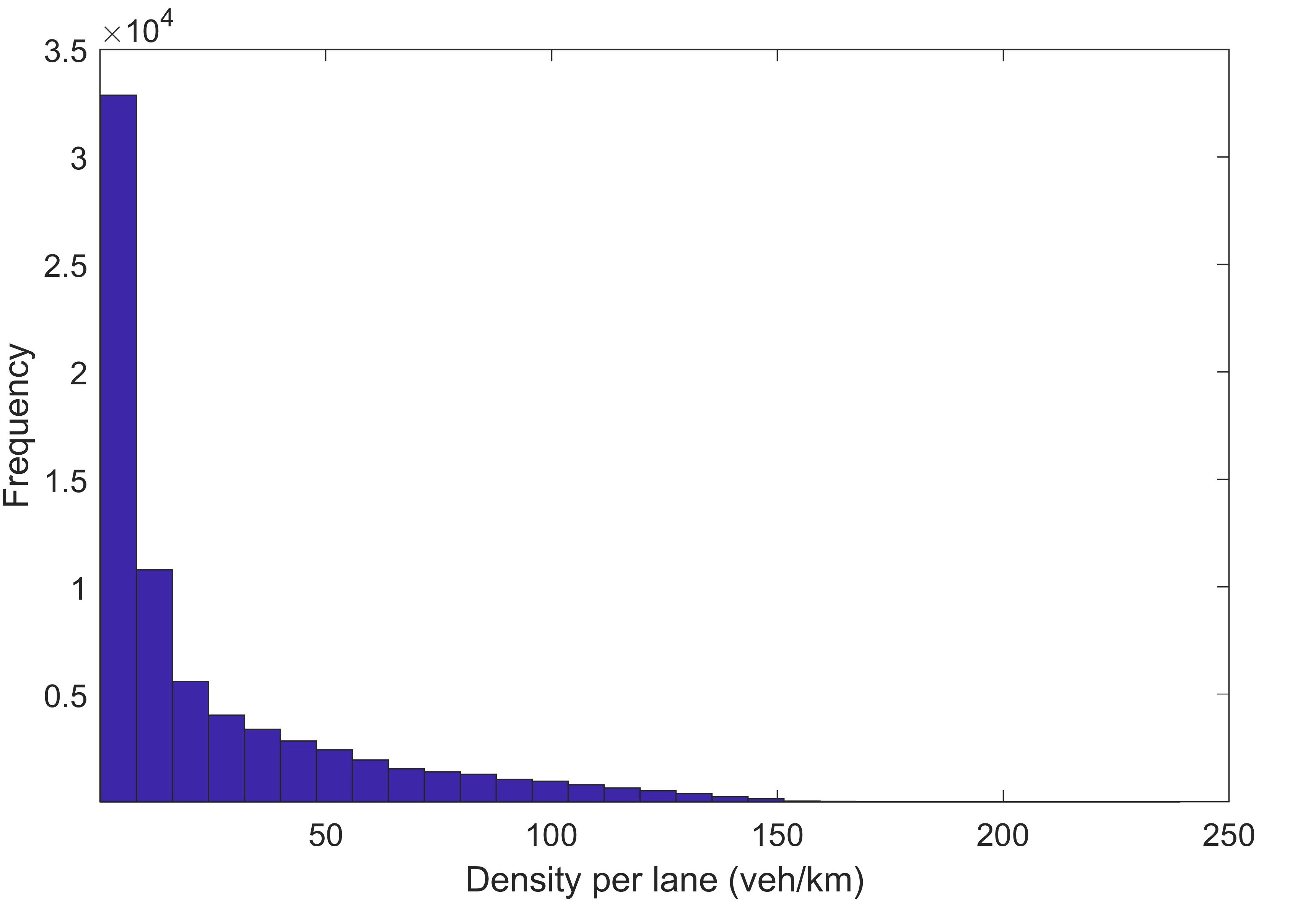}
		\label{histogram density per lane}}
	\hfil
	\subfloat[]{\includegraphics[width=2.5in]{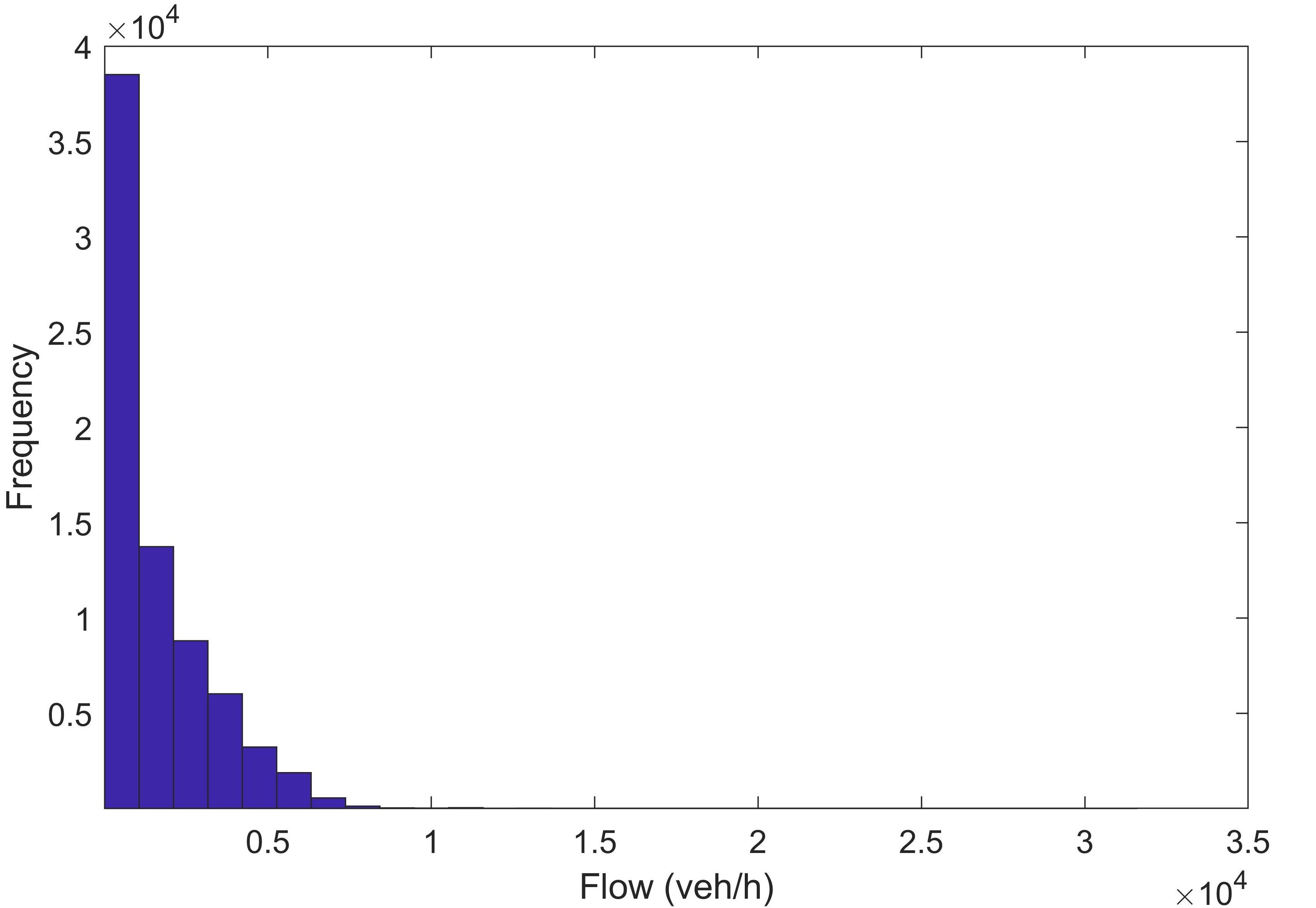}
		\label{histogram total flow}}
	\hfil
	\subfloat[]{\includegraphics[width=2.5in]{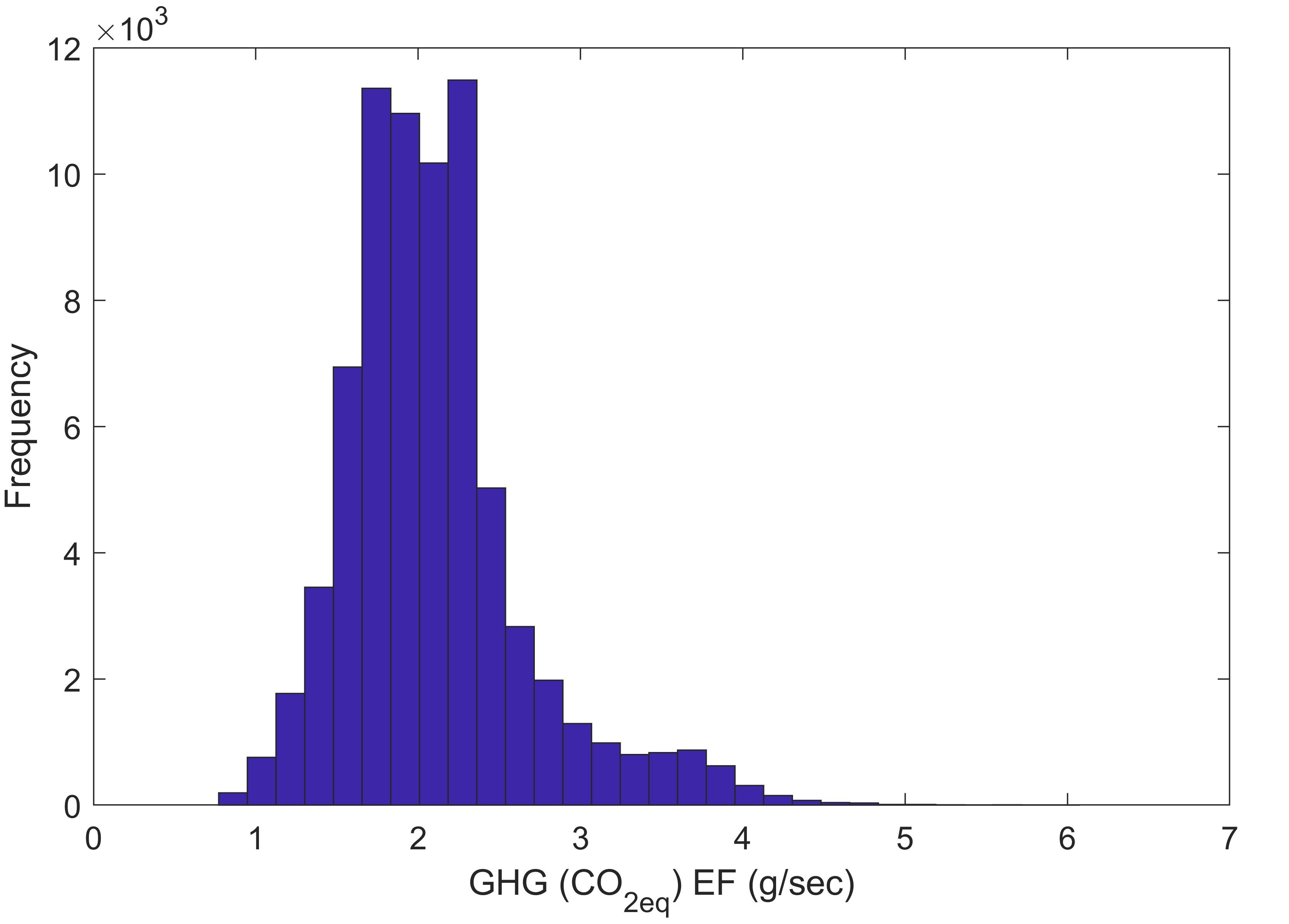}%
		\label{histogram CO2}}
	\hfil
	\caption{Histogram of \textcolor{black}{the traffic and environmental indicators in the network} a) speed, b) density per lane, c) flow, and d) GHG ERs (in $\text{CO}_{2\text{eq}}$ g/sec) \textcolor{black}{\cite{Lama_Prediction2020}}}
	\label{histograms}
\end{figure}

The number of vehicles varies from 2,437 to 6,988, representing 0.7 to 2 times the actual demand in the year 2014. The departure time distributions employed to generate the data are exponential, uniform, and normal. Figure \ref{histograms} illustrates the statistical analysis of the three profound traffic variables, speed, flow, and density in addition to the GHG ER in the data set employed for training the LSTM network. Figure \ref{histogram speed} shows that the mode and average are 40 km/h and 56.16 km/h, respectively. Speed range is from 0 to 80 km/h. Density (veh/km.lane) as in Figure \ref{histogram density per lane} and flow (veh/h) as in {Figure \ref{histogram total flow}} \textcolor{black}{represent different} traffic conditions \textcolor{black}{due to the wide range of their values}. Finally, GHG ER (in $\text{CO}_{2\text{eq}}$) starts from less than 1 g/sec to more than 5 g/sec as in {Figure \ref{histogram CO2}}. 

\subsubsection{Distributed routing framework}
\label{Routing framework description}
\textcolor{black}{A dynamic distributed routing framework, End-to-End Connected Autonomous Vehicles (E2ECAV) \cite{FarooqandDjavadian}, is used for testing the proposed routing strategies. E2ECAV is based on a network of Intelligent Intersections (I2s) that can dynamically route connected and automated vehicles (CAVs). Two types of communication are employed, the vehicle to infrastructure (V2I) and infrastructure to infrastructure (I2I). Via the local communication between the agents, the I2s develop a coherent view of the network. For further description of E2ECAV and various applications we refer the readers to \cite{djavadian2020multi,FarooqandDjavadian,tu2019quantifying,Djavadian2018Distributed,alfaseeh2018impact}.}

\subsection{Routing strategies}
\label{Routing strategies}
Three strategies are examined for both myopic and anticipatory routing, as illustrated in Table \ref{scenarios investigated}. Equation \ref{Optimization process} presents the general formula followed based on the \textcolor{black}{routing objective, TT, GHG, or TT\&GHG}. Where $T_i$ is travel time on link $i$ and $E_i$ is emissions on link $i$ that is the GHG (in $\text{CO}_{2\text{eq}}$) in our case; $n$ is number of links of a path $k$; and $W_t$ and $W_e$ are the weights associated with travel time and emissions, respectively. \textcolor{black}{($T$ and $E$) in \textbf{{Equation \ref{Optimization process}}} are of different units. When multi-objective routing is applied, converting them to a consistent unit (e.g.monetary value) is the solution for a realistic application. The weights in \textbf{{Equation \ref{Optimization process}}} are used for the aforementioned normalization process.} 
The routing strategies considered are illustrated in Table \ref{scenarios investigated}. Every routing strategy is run for five replications of different seeds to account for the stochasticity. When the routing strategy is myopic, current time step $t$ of the variables in concern are adopted. While when the routing strategy is anticipatory, predicted values at $t+1$ of the developed predicted models are employed. For instance, routing strategy TT$_m$ means that TT is obtained based on the current time step $t$ value of speed following Equation \ref{TT_cost}. While TT$\&$GHG$_a$ is the anticipatory routing strategy of when the predicted TT and GHG of time step $t+1$ (of the best trained LSTM networks) are the routing objectives. The TT cost of links at every minute follows Equation \ref{TT_cost}. The link GHG cost in this study is chosen based on the analysis of the different GHG costing approaches as demonstrated in Equation \ref{GHG_cost1}, \ref{GHG_cost2}, \ref{GHG_cost3}, \ref{GHG_cost4}, and \ref{GHG_cost5}.

\begin{equation}
min \displaystyle \Bigg\{\sum_{i=1}^n W_t.T_i+ \sum_{i=1}^n W_e.E_i\Bigg\}
\label{Optimization process}
\end{equation}

\begin{table}[!h]
	\caption{Routing strategies considered}
	\begin{center}
		\small
		\begin{tabular}[!t]{l  c c}
			\hline
			Scenario name & \textcolor{black}{Routing strategy} & Equation \\
			\hline
			\hline
			\centering
			TT$_m$ & Myopic & \ref{TT_cost}\\
			\hline
			GHG$_m$ & Myopic & \ref{GHG_cost5}\\
			\hline
			TT\&GHG$_m$ & Myopic & \ref{GHG_cost5} and \ref{TT_cost} \\
			\hline
			TT$_a$ & Anticipatory & \ref{TT_cost}\\
			\hline
			GHG$_a$ & Anticipatory & \ref{GHG_cost5}\\
			\hline
			TT\&GHG$_a$ & Anticipatory &  \ref{GHG_cost5} and \ref{TT_cost}\\
			\hline
			\hline
		\end{tabular}
	\end{center}
	\label{scenarios investigated}
\end{table}

\section{Results and discussion}
\label{Discussion and results}
Here we compare the results of anticipatory routing strategies with the myopic strategies. To achieve this, predictive models of the related variables i.e. GHG cost and speed, are required. The GHG cost related variable LSTM based predictive model is developed in \cite{Lama_Prediction2020}, but the major findings are shared in the following related sections. The speed predictive model is developed in this work in Section \ref{Development of the predictive models}. The comparison between the routing strategies is conducted in Sections \ref{Routing strategies analysis} to \ref{Network level analysis}.

\subsection{Development of the predictive models}
\label{Development of the predictive models}
Before developing the predictive models, the most representative GHG costing approach is investigated in Section \ref{GHG costing approaches_analysis}. A comprehensive correlation analysis is applied, for the GHG ER and speed, in Section \ref{Correlation analysis of the predicted variables} followed by the predictive models developed in Section \ref{GHG and speed predictive models}.

\subsubsection{GHG costing approaches}
\label{GHG costing approaches_analysis}
\begin{figure}[!h]
	\centering
	\subfloat[]{\includegraphics[width=3.2in]{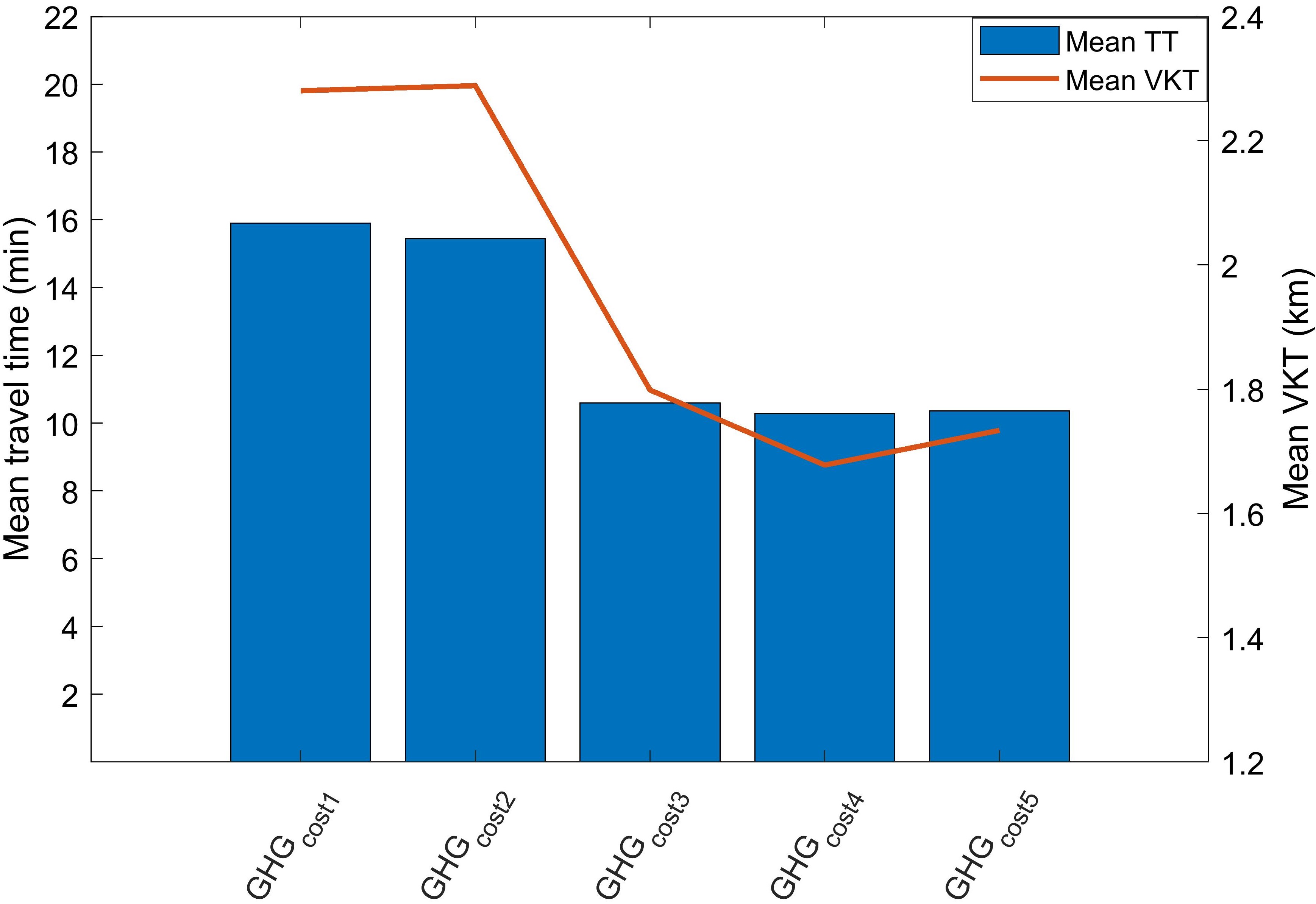}%
		\label{TT_VKT_for_GHG_costing_strategies}}
	\hfil
	\subfloat[]{\includegraphics[width=3.2in]{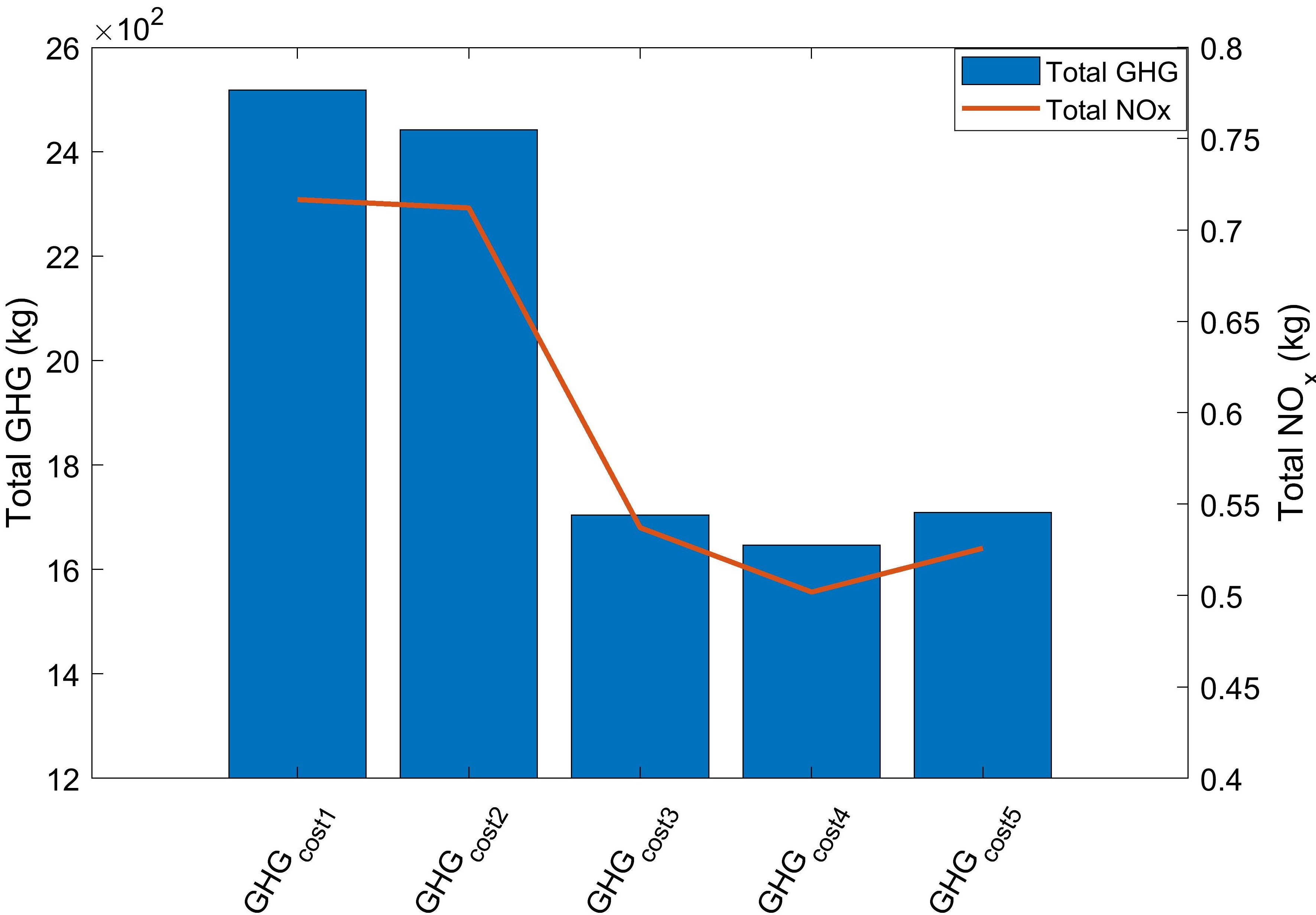}%
		\label{GHG_NOx_for_GHG_costing_strategies}}
	\hfil
	\caption{Performance indicators of the different GHG costing strategies, a) mean TT and VKT and b) total GHG and NOx}
	\label{TT_VKT_GHG_NOx_For_GHG_costing_strategies}
\end{figure}

This analysis is applied for GHG$_m$, \textcolor{black}{as myopic routing is the base case and to illustrate which GHG costing approach is the most suitable for our application.} Figure \ref{TT_VKT_for_GHG_costing_strategies} shows that normalizing based on the number of lanes as in GHG$_{cost2}$ and GHG$_{cost4}$ is associated with a slight enhancement compared to GHG$_{cost1}$ and  GHG$_{cost3}$, respectively. This is due to the different number of lanes of the links in the network that makes the total GHG not reflective of the actual conditions. \textcolor{black}{For the costing approach GHG$_{cost1}$ and GHG$_{cost3}$, if the GHG cost of link $a$ with two lanes and $b$ with four lanes are 70 and 100 grams, respectively, link $a$ would be prioritized. However, dividing by the number of lanes shows that link $b$ should be prioritized based on GHG$_{cost2}$ or GHG$_{cost4}$.} A reduction in the average TT of around 3\% in both GHG$_{cost2}$ and GHG$_{cost4}$ compared to GHG$_{cost1}$ and GHG$_{cost3}$, respectively, is observed. Using the total GHG on links \textcolor{black}{of the GHG$_{cost1}$ costing approach,} triggers the highest average TT, average VKT, total GHG, and total NOx of 16 minute, 2 km, 2,518 kg, and 0.716 kg, respectively. The explanation is related to not considering any of the traffic characteristics on the link, such as speed, density, or flow. \textcolor{black}{100 gram of the total GHG emission, GHG$_{cost1}$, can be for two dramatically different sets of traffic characteristics. The first condition can be for a very congested short link of low capacity, while the other condition can be for an uncongested long link of high capacity.} 

GHG$_{cost3}$ illustrates a significant reduction in the average TT, average VKT total GHG, and total NOx of 33\%, 21\%, 32\%, and 25\%, respectively, when compared to GHG$_{cost1}$. The main justification is the high temporal resolution adopted, compared to GHG$_{cost1}$. Furthermore, giving higher weights to the most recent seconds contributes to the improvements \textcolor{black}{in terms of average TT, average VKT total GHG, and total NOx}. \textcolor{black}{The closer to the prediction update, the more weight is used for the GHG produced by the vehicles.} \textcolor{black}{GHG$_{cost5}$, which is the marginal cost of 1 vehicle based on the GHG ER and TT on the studied link, is very much comparable to the  GHG$_{cost4}$ in terms of the performance. It is despite the lower temporal resolution utilized by GHG$_{cost5}$ (1 min) when compared to  GHG$_{cost4}$ (1 second)}. GHG$_{cost5}$ is the most representative and suitable costing approach as it is based on the reflective ER and speed on the studied link as illustrated in Equation \ref{GHG_cost5}. Due to the quasi-convex relationship between the speed and GHG ER \textcolor{black}{\cite{djavadian2020multi}}, too low or too high speed will trigger higher emission rates. \textcolor{black}{Links with too high speed are associated with high GHG ERs, but less travel time. Links with too low speed contribute adversely not only to the GHG ER, but also to the travel time on the links.} It is expected that the GHG marginal cost will search for the optimal combination of speed, VKT, and GHG ER to satisfy the optimization Equation \ref{Optimization process}. GHG$_{cost5}$ is the cost used for the myopic and anticipatory routing strategies in this study.

\subsubsection{Correlation analysis of the predicted variables}
\label{Correlation analysis of the predicted variables}
\begin{figure*}[ht]
	\centering
	\includegraphics[width=5in]{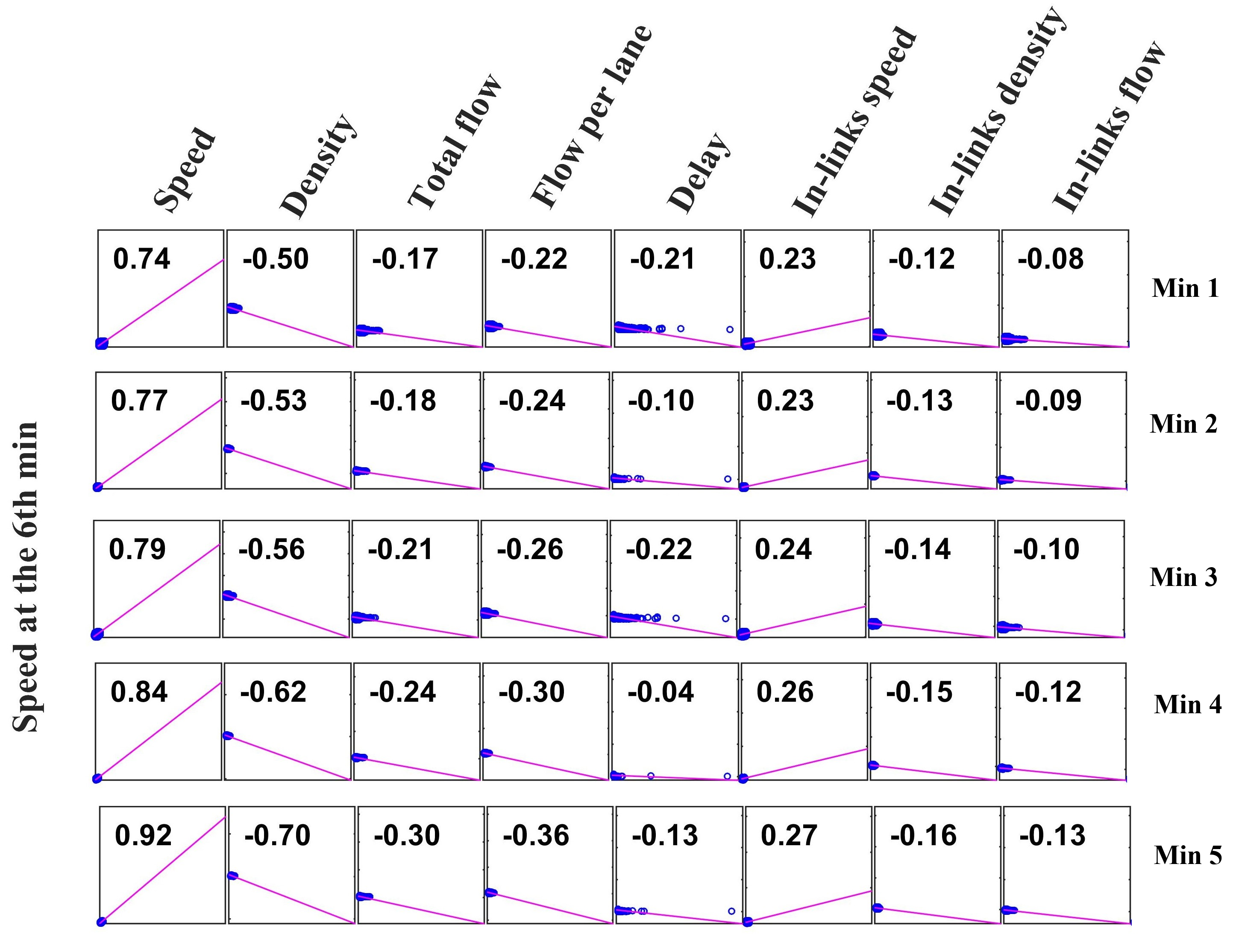}
	\caption{Correlation between the \textcolor{black}{traffic and environmental} variables at \textcolor{black}{the previous five minutes (top (1) to bottom (5)} with speed at the sixth minute}
	\label{Correlation_analysis_speed}
\end{figure*}

This analysis is \textcolor{black}{an essential} step for reliably selecting the predictors and number of sequences of the LSTM models developed. For the GHG ER correlation analysis, details can be found in \cite{Lama_Prediction2020}. \textcolor{black}{A comprehensive list of variables has been assessed to develop reliable predictive models. Not only the link characteristics (speed, density, flow, and delay (difference between free flow travel time and actual travel time)), but also the in-links characteristics (speed, density, and flow) are included for the correlation analysis.} 
The in-links characteristics are included to enhance the spatial dimension. \textcolor{black}{The traffic state at time $i$ on the up steam links implies how the traffic condition will be on the downstream links at time $i+1$.}
\textcolor{black}{Five time steps (minutes) are considered for this analysis as within the aforementioned period modifications in the traffic state is detected.}
\textcolor{black}{The maximum link length in the network is around 450 meters. \textcolor{black}{The speed range is from} 0 to 80 km/h. Under the free flow traffic condition, the maximum travel time required to traverse a link is around 0.8 minute.} \textcolor{black}{The GHG ER at the sixth minute is highly correlated with speed, GHG ER, density, and in-links speed over the last five minutes, followed by the rest of the variables, as shown in \cite{Lama_Prediction2020}. The GHG emission estimation \cite{MOVES} depends primarily on speed, which justifies the strong correlation. Speed and density have a more or less monotonically decreasing relationship \cite{papacostas1993transportation}, which explains the strong correlation between the GHG ER and density.}
In terms of the correlation analysis for speed, Figure \ref{Correlation_analysis_speed} shows the linear correlation between speed at the sixth minute and both the traffic and environmental indicators of the previous five minutes. Not only the indicators on the studied links are considered, but also the in-links characteristics as in the case of the GHG ER correlation analysis to better reflect on the spatial dimension. Figure \ref{Correlation_analysis_speed} shows an increase in the correlation factor of all the variables and minutes, except for delay. The top four highly correlated variables with speed at the sixth minute are speed, density, flow per lane, and in-links speed. The high correlation coefficient with density of 0.70 and 0.36 with flow is based on the correlation between the three variables shown in the transportation fundamental diagrams, speed, density, and flow \cite{papacostas1993transportation}. Speed and density are associated with a monotonically decreasing relationship \cite{papacostas1993transportation}.

\subsubsection{GHG and speed predictive models}
\label{GHG and speed predictive models}
\begin{figure}[!h]
	\centering
	\subfloat[]{\includegraphics[width=3.2in]{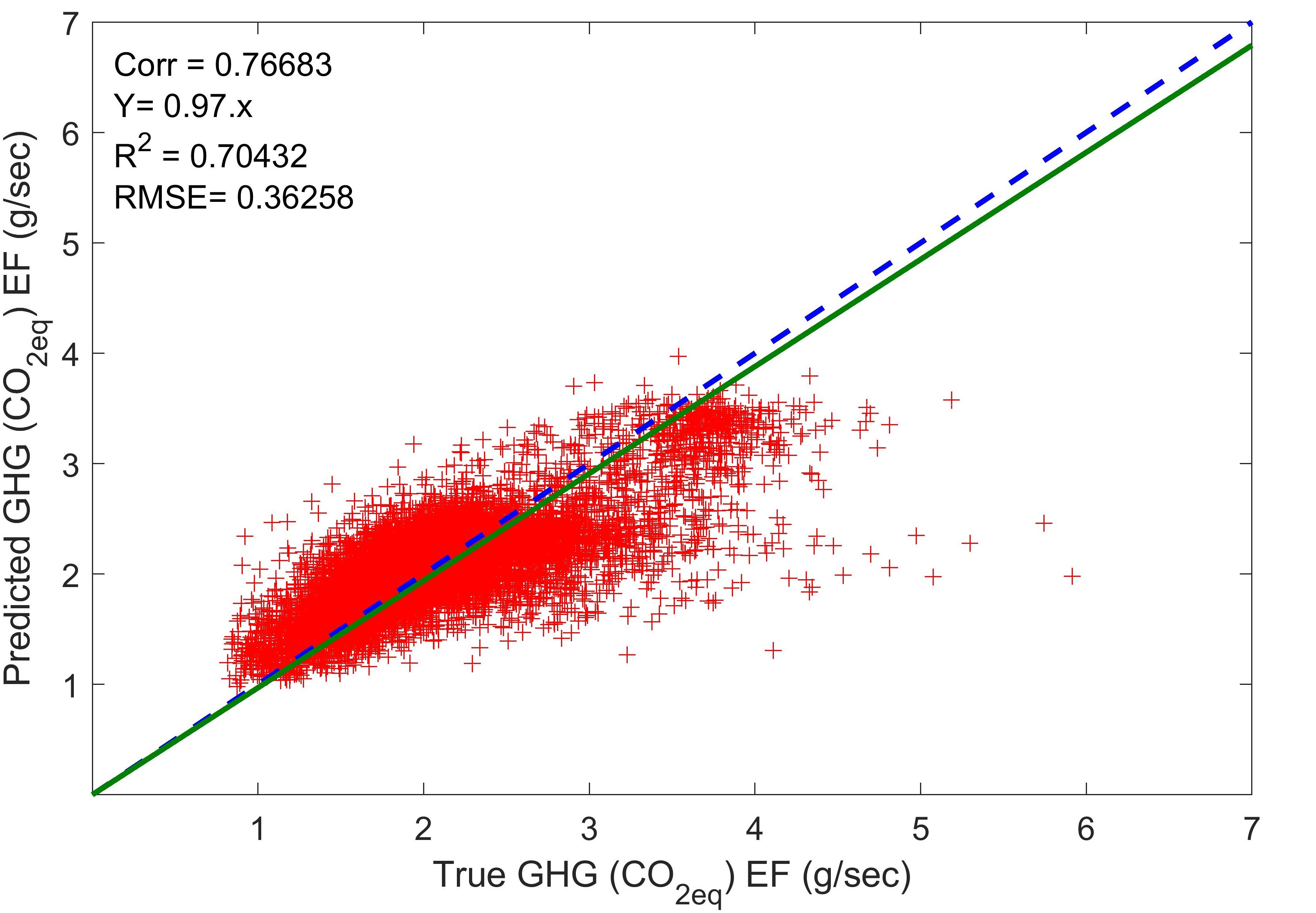}%
		\label{GHG_best_LSTM_model}}
	\hfil
	\subfloat[]{\includegraphics[width=3.2in]{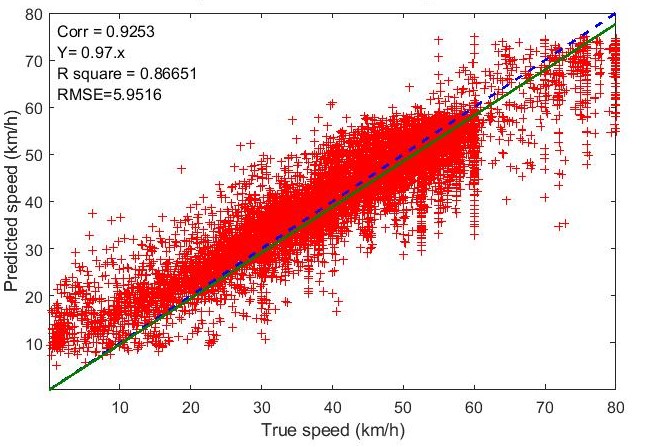}%
		\label{Speed_best_LSTM_model}}
	\hfil
	\caption{True vs. predicted of the best LSTM models for GHG ERs (in $\text{CO}_{2\text{eq}}$ g/sec) \cite{Lama_Prediction2020} and speed, respectively. a) GHG ERs and b) speed}
	\label{GHG_speed_final_LSTM_models}
\end{figure}

For each of the predicted variables, GHG ER and speed, a comprehensive list of predictors and number of sequences have been examined. The best predictive LSTM networks of both the GHG ER and speed are of two hidden layers while the hyper-parameters are systematically tuned. In terms of the predictors, for GHG ER forecasting \cite{Lama_Prediction2020} the best set is of three sequences of speed, GHG ER, density, and in-links speed. For speed prediction, the best setting is associated also with the three sequences of speed, density, and in-links speed. In terms of the hyper-parameters, two solvers are considered, the adaptive moment estimation (Adam) methods \cite{kingma2014adam} and the stochastic gradient descent with momentum (sgdm) \cite{robert2014machine}. \textcolor{black}{Several hyper-parameters are considered for tuning. The initial learning rate, max epochs, learning rate drop factor, momentum, learning rate drop period, number of hidden units of the first LSTM (hidden) layer, and the number of hidden units of the second (LSTM) layer when used are the hyper-parameters tuned.} 
The first stage of tuning is manual, which is adopted to narrow the search range of the optimal hyper-parameters. The next stage is systematically using the Bayesian optimization \cite{wu2019hyperparameter, Lama_Prediction2020}, which employs a narrowed search range obtained from the manual tuning stage. The training results of the best LSTM predictive networks of the GHG ER and speed are shown in Figure \ref{GHG_speed_final_LSTM_models}, respectively.

\textcolor{black}{For} the prediction performance, \textcolor{black}{four indicators are deployed.} The correlation coefficient between observed and predicted GHG ERs (in $\text{CO}_{2\text{eq}}$ g/sec), the fit to the ideal straight curve reflecting on the precision, R$^2$ statistics, and the RMSE reflecting on the accuracy are used in this study \cite{Lama_Prediction2020}. The correlation coefficient of GHG ER and speed prediction is 0.77 and 0.92, respectively. The RMSE of the GHG ER and speed predictive models is 0.36 gram and 5.95 km\/h, respectively. The performance of the speed predictive model is noticeably better than the one of the GHG ER. This is due to the complicated relationship between the GHG ER and the predictors, while speed has more straight forward relationships with the predictors used. \textcolor{black}{For instance, speed and density have a monotonically decreasing relationship \cite{papacostas1993transportation}, while the GHG ER has a quasi-convex relationship with the most important predictor (speed) \cite{djavadian2020multi}.}

It is noticed in both Figure \ref{GHG_best_LSTM_model} and \ref{Speed_best_LSTM_model} that true values higher than 4 g/sec and 60 km/h, respectively, are not predicted with a high level of accuracy. This stems from the fact that the \textcolor{black}{frequency of the data points reflecting these conditions is much less compared to the other conditions.} as in Figure \ref{histograms}. \textcolor{black}{The data points reflecting on the GHG ER greater than 4 g/sec, as in Figure \ref{histogram CO2}, is only 0.007\% of the total GHG ER data points. The high GHG is also associated with either low or high speed. \textcolor{black}{It is probably because of the quasi-convex relationship of GHG ER with speed} \cite{djavadian2020multi}}. Similarly, the data points of speed higher than 60 km/h and less than 20 km/h \textcolor{black}{represent only 0.11\% of the total data points, while the data points of speed between 20 and 60 km/h represent 0.89\% as in Figure \ref{histogram speed}.}

\subsection{Routing strategies analysis}
\label{Routing strategies analysis}
\begin{figure}[!h]
	\centering
	\subfloat[]{\includegraphics[width=3.2in]{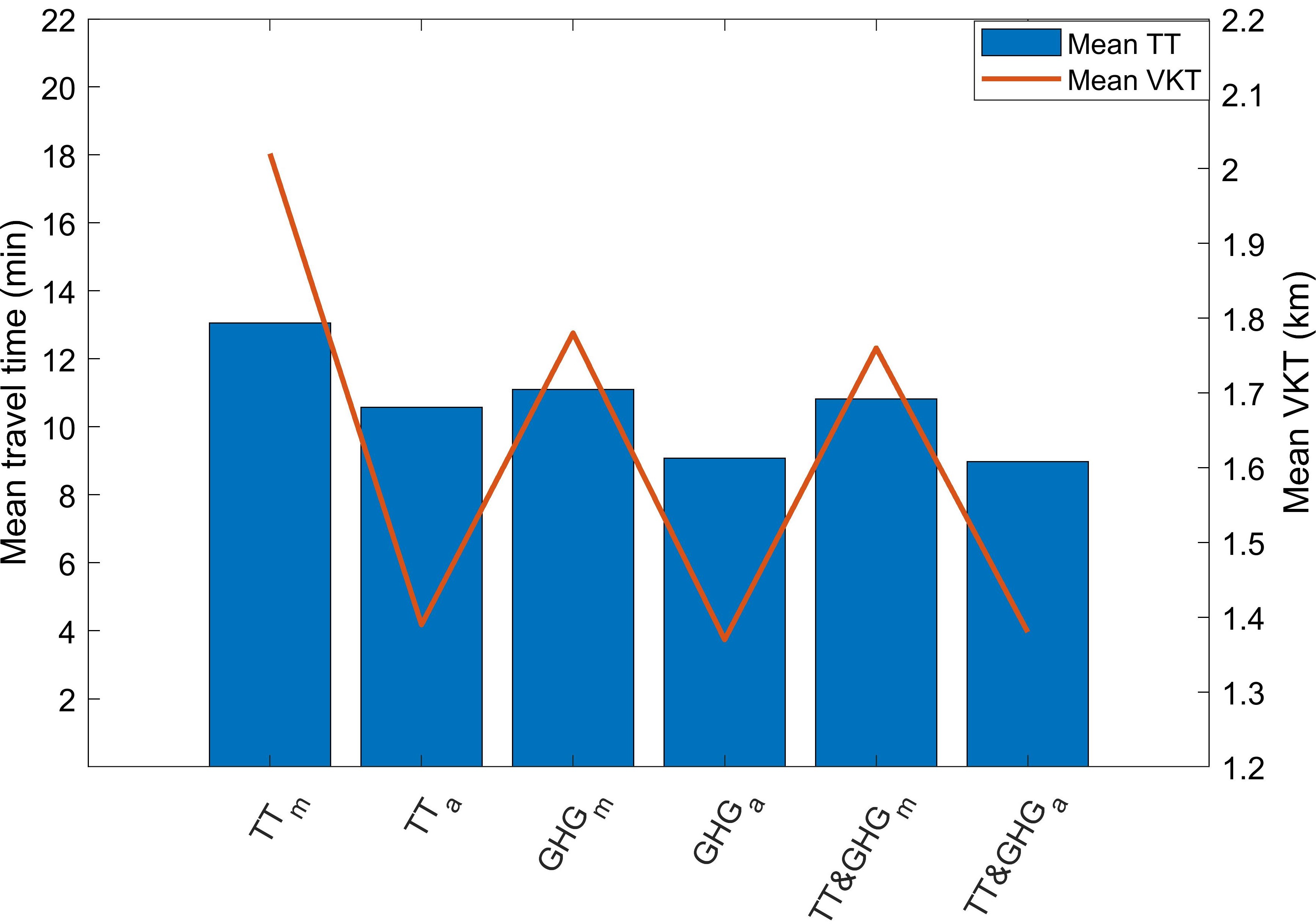}%
		\label{TT_VKT_routing_strategies}}
	\hfil
	\subfloat[]{\includegraphics[width=3.2in]{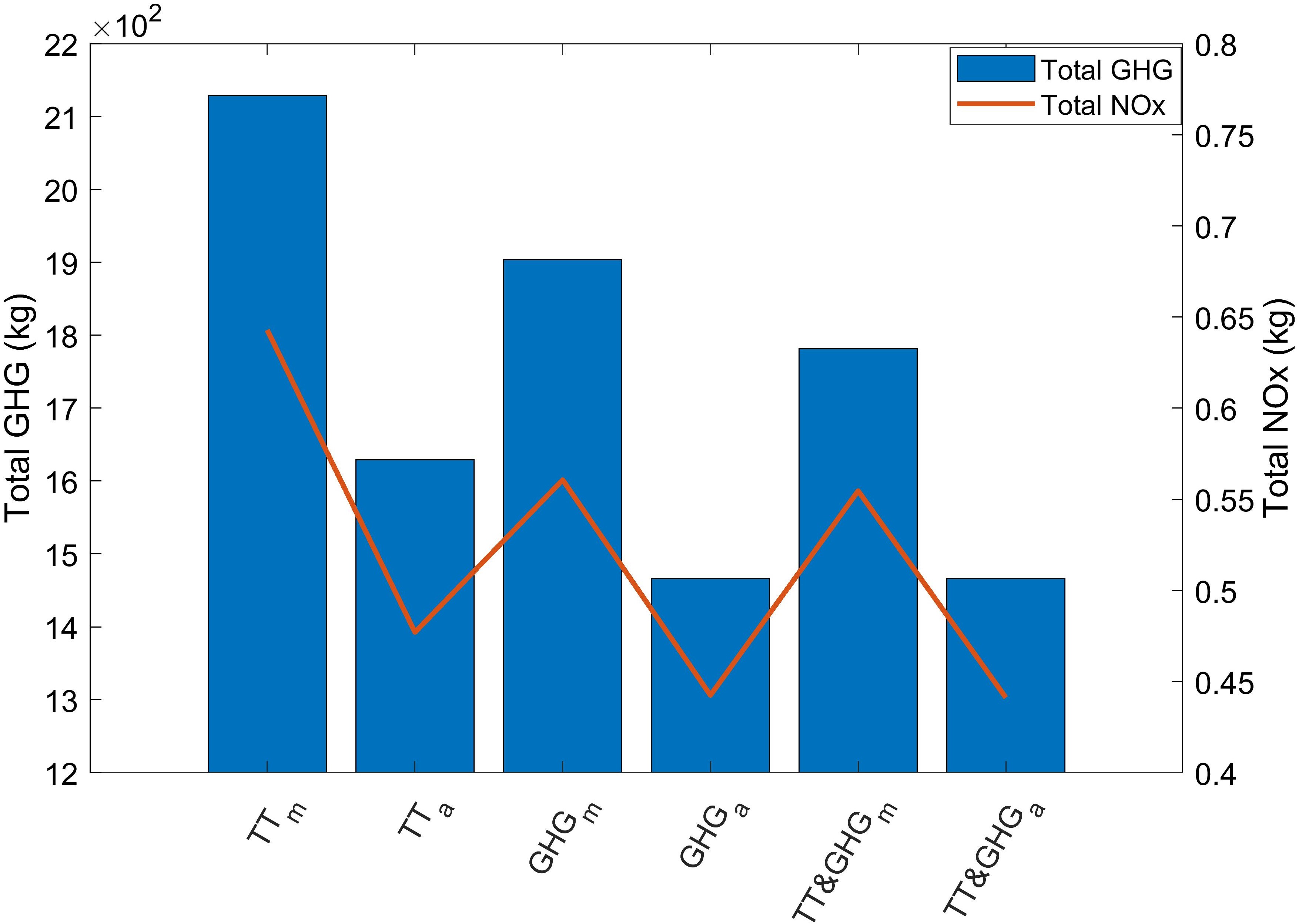}%
		\label{GHG_NOx_routing_strategies}}
	\hfil
	\caption{Performance indicators of the different routing strategies, a) mean TT and VKT and b) total GHG and NOx}
	\label{Mean_TT_mean_VKT_total_GHG_NOx_routing_strategies}
\end{figure}

As illustrated in Table \ref{scenarios investigated}, six routing strategies are analyzed using the E2ECAV \cite{FarooqandDjavadian} routing framework. Single and multi-objectives are considered while routing is myopic and anticipatory. Mean TT, mean VKT, total GHG, and total NOx are the performance indicators taken into account and the results are shown in Figure \ref{Mean_TT_mean_VKT_total_GHG_NOx_routing_strategies}. NOx is the pollutant reflecting on the public health \cite{alfaseeh2020multifactor}. It is \textcolor{black}{of high importance} to include the NOx as a performance factor to assess the impact of the different routing strategies. It is \textcolor{black}{important} to note that \textcolor{black}{logical} constraints have been included for a realistic application while estimating the cost on links. When either the predicted GHG ER or speed of time step $t+1$ is negative, the value is set to zero. Predicted speed of time step $t+1$ that surpasses the link speed limit is set to the link speed limit. With regards to the results, Figure \ref{Mean_TT_mean_VKT_total_GHG_NOx_routing_strategies} illustrates that the anticipatory routing strategies outperform the myopic ones whatever the routing objective is. \textcolor{black}{The justification of this outcome is threefold. Firstly, the reflective costing approach of the GHG emissions (when GHG is part of the optimization process), which optimizes not only the GHG ER, but also the travel time implicitly to avoid re-routing is used. The marginal cost prioritizes the links of speed close to the optimal value based on the quasi-convex relationship between the GHG ER and speed \cite{djavadian2020multi}.
	Secondly, the sophisticated predictive models developed based on high resolution data points are adopted. Lastly, \textcolor{black}{by taking into account the traffic conditions and their evolution, routing is more proactive than simply being reactive to the current conditions}.}


The performance trend for myopic routing strategies is similar to that of the anticipatory routing strategies. From the worst to the best, TT is followed by GHG and TT\&GHG in terms of the four performance indicators. Whether it is myopic or anticipatory routing, when TT is \textcolor{black}{optimized} the worst performance indicators are observed compared to when \textcolor{black}{the routing objective is GHG or TT\&GHG}. The justification is that when TT is the objective all that matters is the TT spent on the links regardless of the VKT, GHG and NOx indicators. The vehicles are distributed in the network to achieve the least TT, but this comes at the cost of longer distances travelled and more GHG and NOx produced. Nevertheless, GHG$_m$ and GHG$_a$ introduce a decrease in the total GHG produced of 11\% and 6\% compared to TT$_m$ and TT$_a$, respectively. This improvement can be directly linked to the marginal costing approach as in Equation \ref{GHG_cost5}, which takes into account not only the GHG ER, but also the speed on links. In other words when the routing objective is GHG, the chosen links are defined based on the best combination of the GHG ER and speed that minimizes the cost following Equation \ref{Optimization process}. The relationship between GHG and speed is quasi-convex \textcolor{black}{\cite{djavadian2020multi}} in which too low or too high speed will contribute to a higher GHG ER and higher GHG cost eventually. The multi-objective TT\&GHG$_m$ outperforms both TT$_m$ and GHG$_m$ in terms of the whole performance indicators. A reduction in average TT, average VKT, total GHG, and total NOx of 17\%, 13\%, 16\%, and 14\%, respectively, is observed for the TT\& GHG$_m$ strategy compared to TT$_m$ strategy. The reduction in the performance indicators of TT\&GHG$_m$ is marginal compared to the GHG$_m$ routing strategy. TT\&GHG routing objective controls the TT cost and does not allow it to neglect the GHG objective. Thus, TT is reduced as long as it still satisfies the objective of reducing GHG. Paths of longer distances are not chosen as this triggers the increase in both the GHG and NOx produced like in TT routing strategy. Comparing the best anticipatory routing strategy, TT\&GHG$_a$ to the TT$_m$, a reduction in average TT, average VKT, total GHG, and total NOx of 17\%, 22\%, 18\%, and 20\%, respectively is noticed. With reference to the NOx variable, it has been found that the relationship between NOx and speed is quasi-convex \cite{djavadian2020multi}. Moreover, \textcolor{black}{\textcolor{black}{previous studies have confirmed that at high speeds NOx is sensitive to aggressive driving} \cite{tu2019quantifying}.} Figure \ref{GHG_NOx_routing_strategies} shows that when GHG is part of the routing objective (GHG or TT\&GHG), the NOx produced is less compared to when TT is the objective regardless of the routing protocol, myopic and anticipatory. It can be concluded that the additional time and longer trips experienced by vehicles in the case of myopic and anticipatory routing contribute to the increase in the GHG and NOx emissions produced.

\subsection{Path analysis}
\label{Path analysis}
\begin{figure}[!h]
	\centering
\includegraphics[width=6in]{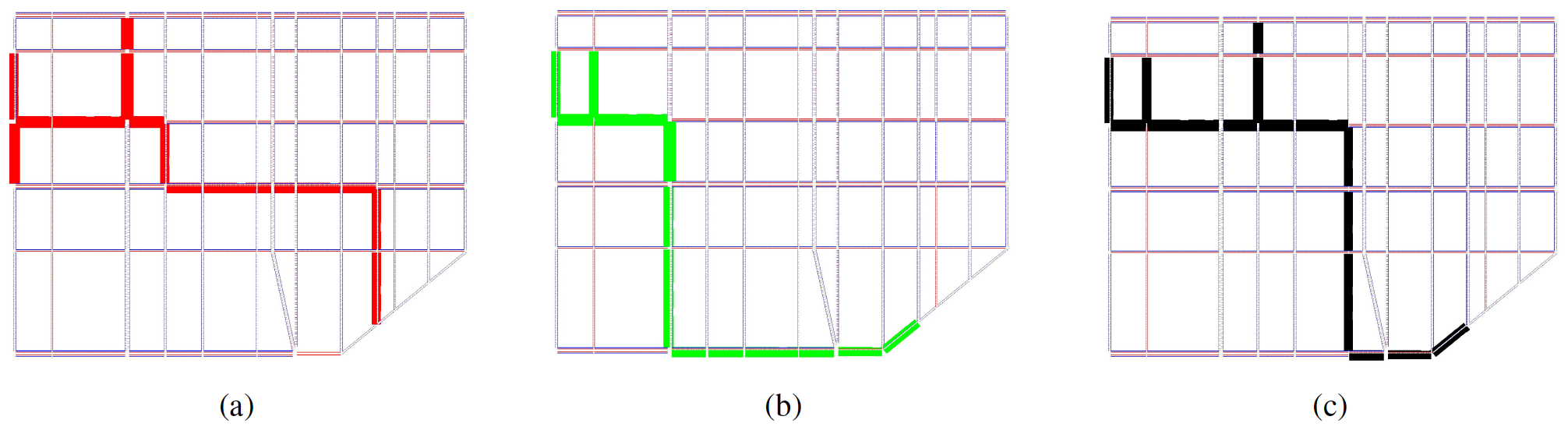}%
	\caption{Myopic routes of a random vehicle of different routing objectives. a) TT, b) GHG, and c) TT\&GHG}
	\label{veh1660_routes_myopic}
\end{figure}

\begin{figure}[!h]
	\centering
\includegraphics[width=6in]{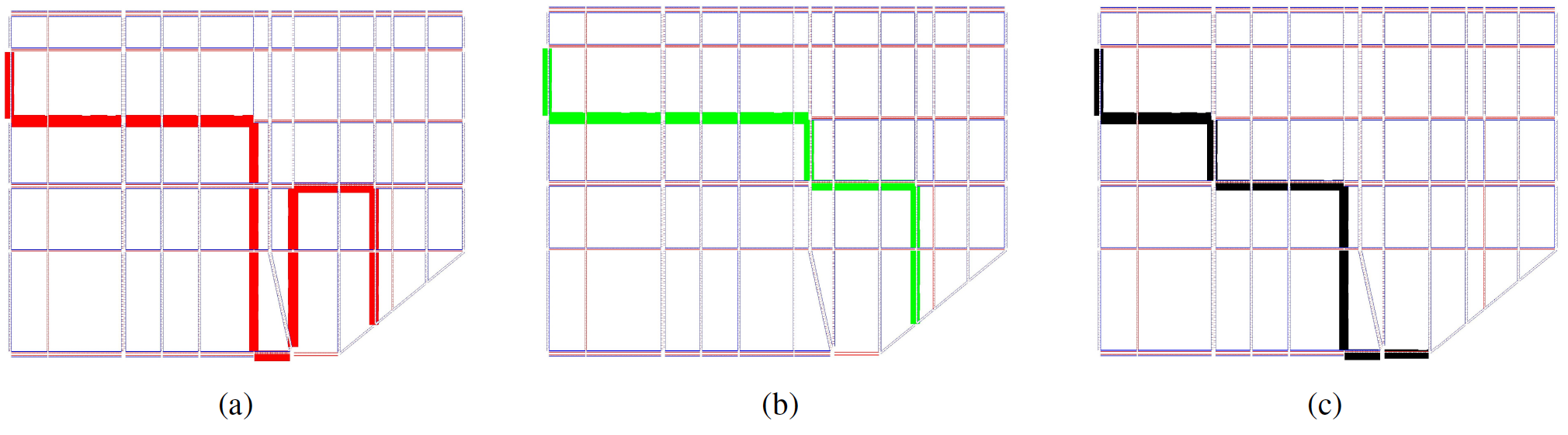}%
	\caption{Anticipatory routes of a random vehicle of different routing objectives. a) TT, b) GHG, and c) TT\&GHG}
	\label{veh1660_routes_anticipatory}
\end{figure}

For this analysis, one vehicle is chosen randomly and its myopic and anticipatory paths with different routing objectives are investigated. Comparing the myopic from Figure \ref{veh1660_routes_myopic} with the anticipatory routes in Figure \ref{veh1660_routes_anticipatory} for each of the objectives, shows that there is more re-routing in the former contributing to longer trips and probably more time in the network as illustrated in Section \ref{Routing strategies analysis}. This \textcolor{black}{main explanation is} that the cost of links is based on the current traffic conditions and does not consider the evolution of traffic in future. This analysis supports the findings in Figure \ref{TT_VKT_routing_strategies}, which demonstrates the decrease in average TT and VKT while anticipatory routing is adopted. Whether the routing protocol is myopic or anticipatory, comparing the length of the path of the TT to GHG and TT\&GHG routing strategies illustrates that the length of the path of TT strategy is the longest. The main justification is that when TT is the objective, vehicles are distributed in the network utilizing uncongested links to achieve the least TT regardless of the distance travelled, the GHG, and NOx produced. Figure \ref{veh1660_routes_anticipatory} illustrates that instead of taking route $a$ (west-east), of two link, route $b$ (south-north), of five links, is chosen. \textcolor{black}{The vehicle traversed an additional distance of around 700 meters when taking route $b$ compared to GHG$_a$ and TT\&GHG$_a$ routing strategies. The speed of route $b$ is around 38 km/h, while the speed of route $a$ is around 1 km/h.} This asserts that vehicles are distributed in the network to links of high speed to minimize the TT regardless of the distance travelled. When TT is the objective, the time spent on the links is optimized, while when GHG is part of routing objective the links of optimal speed are prioritized as long as the GHG marginal cost is minimized. The length of routes of the GHG and TT\&GHG routing strategies is comparable as illustrated in Figure \ref{veh1660_routes_anticipatory}. 

\subsection{Network level analysis}
\label{Network level analysis}
\begin{figure}[!h]
	\centering
	\includegraphics[width=6in]{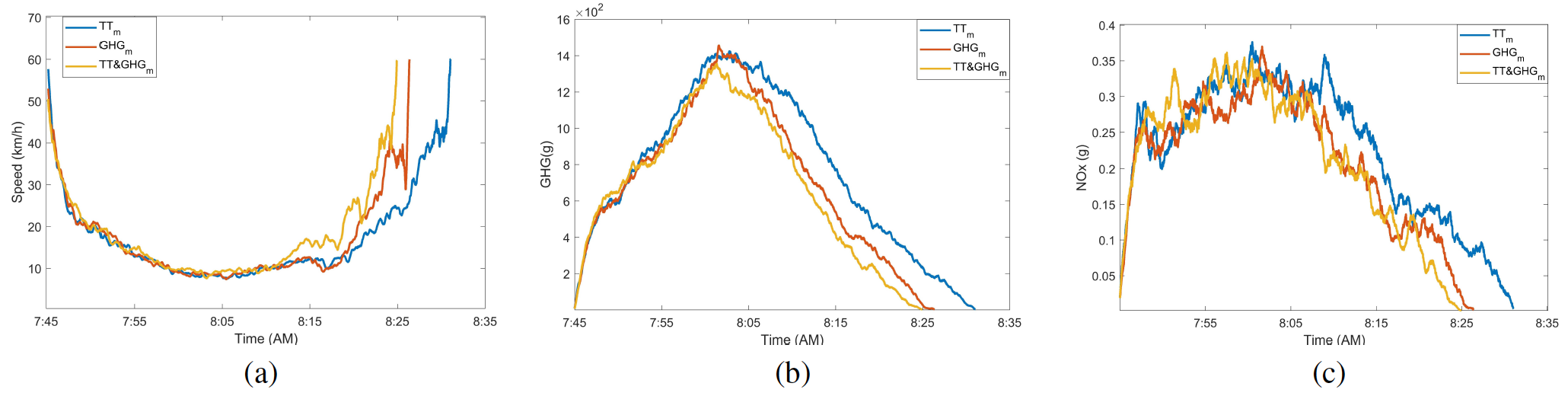}
	\caption{\textcolor{black}{Network level analysis of the myopic routing of different objectives}. a) Average speed, b) GHG, and c) NOx}
	\label{Network_NOx_myopic}
\end{figure}

\begin{figure}[!h]
	\centering
	\centering
	\includegraphics[width=6in]{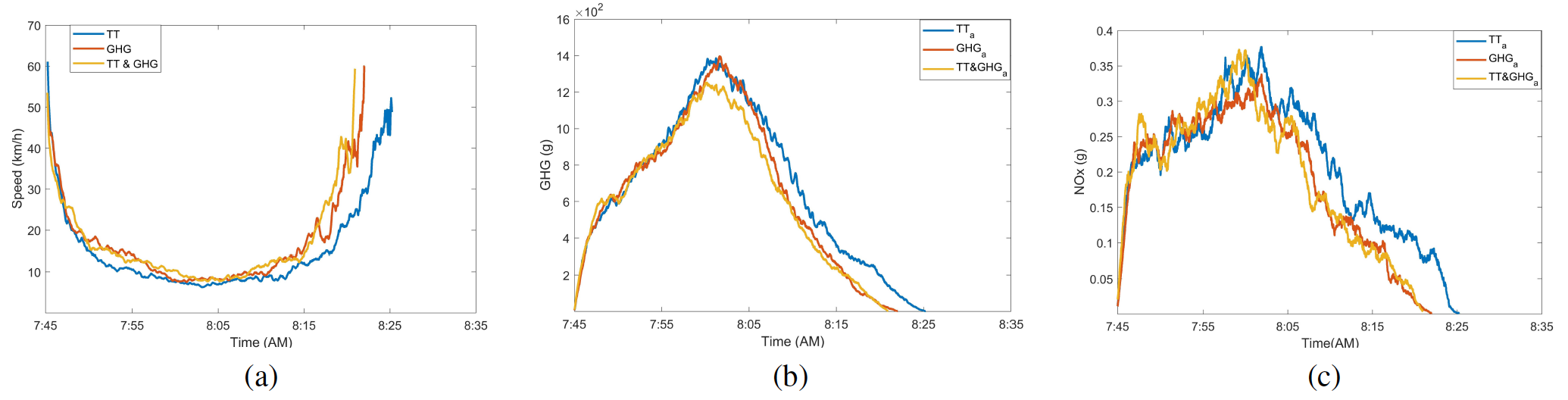}
	\caption{\textcolor{black}{Network level analysis of the anticipatory routing of different objectives}. a) Average speed, b) GHG, and c) NOx}
	\label{Network_NOx_non_myopic}
\end{figure}

To \textcolor{black}{examine the effect} at the network level of the myopic and anticipatory routing strategies while adopting different routing objectives, average speed, GHG, and NOx produced over time are examined. \textcolor{black}{The demand is loaded at 7:45am and the total demand is the network at around 8:00am, which represents the peak of the congestion.} Comparing Figure \ref{Network_NOx_myopic} to Figure \ref{Network_NOx_non_myopic}, shows that the network has been loaded and unloaded quicker in \textcolor{black}{the case of} anticipatory routing. Particularly, for TT$_m$ the vehicles spend 15\% more time in the network compared to TT$_a$. This \textcolor{black}{finding is aligned with one of the TT, VKT,} GHG, and NOx in Figure \ref{TT_VKT_routing_strategies} for TT$_m$ \textcolor{black}{compared to TT$_a$}. Figures \ref{veh1660_routes_myopic} and \ref{veh1660_routes_anticipatory} demonstrate the re-routing in the case of myopic and anticipatory routing strategies, respectively. \textcolor{black}{The additional re-routing noticed in the case of TT$_m$, which triggers longer trip lengths compared to TT$_a$, contributes to the additional time spent in the network as well.} \textcolor{black}{The throughput of TT$_m$ is less compared to GHG$_m$ and TT\&GHG$_m$. It takes around 8\% and 10.6\% less time to load and unload the network for GHG$_m$ and TT\&GHG$_m$, respectively, compared to TT$_m$. The main explanation is that when travel time minimization is the objective, all link options are analyzed and the final cost of travel time does not take into account the VKT as long as the objective is minimized. However, when GHG is the part of the optimization process, the links of optimal speed are prioritized and the re-routing is averted. The GHG marginal cost as in Equation \ref{GHG_cost5} makes sure that vehicles spend the least time and travel the least distance while the objective is minimized which is observed in Figure \ref{veh1660_routes_myopic} and \ref{veh1660_routes_anticipatory}. The average speed till 8:10am of the myopic routing strategies in Figure \ref{Network_NOx_non_myopic} is almost identical for the three routing strategies, while for the anticipatory routing GHG$_a$ and TT\&GHG$_a$ are associated with a slightly higher speed than TT$_a$ as in Figure \ref{Network_NOx_non_myopic}. This asserts the importance and positive impact of the anticipatory routing, which takes into account the future state of traffic conditions in the network. After 8:10am and till the end of the simulation, the increase in speed for GHG$_a$ and TT\&GHG$_a$ compared to TT$_a$ is higher than in the case of the GHG$_m$ and TT\&GHG$_m$ compared to TT$_m$. The reason is that the vehicles reach their destination faster in the case of anticipatory routing compared to myopic routing, which means less vehicles are in the network in the former case.} 

The main difference between Figure \ref{Network_NOx_myopic} and \ref{Network_NOx_non_myopic} is that the period of time vehicles produce GHG is longer in the former as the anticipatory routing includes the future state of the traffic conditions \textcolor{black}{and deals with the changes proactively compared to the myopic routing. In addition, the GHG costing approach takes into account not only the GHG ER, but also the speed and VKT implicitly.} \textcolor{black}{The additional time and VKT experienced by the vehicles in the case of TT$_m$, as shown in Figure \ref{TT_VKT_routing_strategies}, contribute to the higher levels of GHG compared to GHG$_m$ and TT\&GHG$_m$ in Figure \ref{Network_NOx_myopic}, especially after the congestion peak at 8:05am. The number of vehicles in the network is an essential factor to keep in mind. On the other hand, comparing the GHG over time of TT$_a$ to GHG$_a$ and TT\&GHG$_a$, illustrates less variation.} 

NOx over time follows the same trend of GHG over time and is associated with less emissions over time while the anticipatory routing is utilized as in Figure \ref{Network_NOx_non_myopic} compared to Figure \ref{Network_NOx_myopic}. For the NOx over time analysis, as shown in Figures \ref{Network_NOx_myopic} and \ref{Network_NOx_non_myopic}, the high values till around 8:00am are due to the high speed of the uncongested network. After 8:00am the number of vehicles and speed level control the NOx produced. \textcolor{black}{At 8:00am, the complete demand is in the network and the speed is the variable with direct impact on the NOx produced.} It is observed that for TT$_m$ in Figure \ref{Network_NOx_myopic} and TT$_a$ in Figure \ref{Network_NOx_non_myopic}, NOx is higher than the cases when GHG or TT\&GHG is the routing objective. This is because vehicles were directed to longer paths, but of higher speeds to minimize the travel time. NOx is sensitive to aggressive driving \cite{tu2019quantifying}, which makes higher speed links unfavourable. However, the reduction in the TT and VKT experienced in the network in the case of anticipatory routing \textcolor{black}{means higher throughput over time. The higher throughput contributes to less NOx over time.}

\textcolor{black}{It can be concluded that using predicted link cost is associated with significant improvements at the network level. Furthermore, utilizing the GHG marginal cost, which takes into account not only the GHG ER, but also the speed and VKT, is very effective and outperforms the travel time cost in terms of the whole performance indicators.}

\section{Conclusion and potential directions}
\label{Conclusion}
\textcolor{black}{Current eco-routing studies are predominantly associated  with limitations related to the \textcolor{black}{aggregation level} used in the traffic flow models, scale of the case study, centralized routing, and the number of objectives optimized \textcolor{black}{simultaneously} \cite{alfaseeh2020multifactor}.  \cite{djavadian2020multi} and \cite{alfaseeh2019multi} overcame these limitations and applied myopic multi-objective eco-routing strategies in a distributed routing framework with favourable outcomes. However, the technological advancements related to ICT and CAVs have not yet been exploited completely. Hence, this study suggested anticipatory multi-objective eco-routing strategies using a distributed routing system for connected and automated vehicles i.e. E2ECAV \cite{FarooqandDjavadian}. Predictive models of GHG ER and speed were developed and used. A deep learning based time series model i.e. LSTM is trained while systematically tuned. \textcolor{black}{For sequential data, LSTM is known to be} the most powerful recurrent NN architectures \cite{lipton2015critical}. Furthermore, the LSTM model employed here outperformed the commonly used statistical time series model e.g. ARIMA and clustering \cite{Lama_Prediction2020}.} 

The major findings of this study are as follows. Anticipatory routing strategies outperform the myopic ones due to the inclusion of future traffic conditions in the route calculations. The paths of myopic routing strategies demonstrate a high degree of re-routing as the cost does not consider the future traffic  and environmental conditions. Routing based on GHG as the objective is associated with noticeable reduction in average TT, average VKT, total GHG, and total NOx compared to the case where TT is the objective. This stemms from the costing approach for the marginal cost in Equation \ref{GHG_cost5}, that results in the best combination of speed and GHG ER on links that minimizes the GHG cost. \textcolor{black}{Re-routing is minimized when GHG is part of the optimization process as the increase in the VKT has a negative impact on the GHG produced.} The GHG routing objective contributes to less TT, VKT that leads to less GHG, NOx in the network. For myopic and anticipatory routing, TT\&GHG routing strategy introduces a slight enhancement in terms of the four performance indicators compared to GHG routing strategy. Comparison of the former to the latter would be like comparing the system optimal to the User Equilibrium (UE) \cite{papacostas1993transportation}. Comparing the best anticipatory routing strategy, which optimizes not only TT, but also GHG, a reduction in average TT, average VKT, total GHG, and total NOx of 17\%, 21\%, 18\%, and 20\%, respectively, is noticed when compared with the myopic routing aiming at TT minimization. 

\textcolor{black}{For future work, utilizing real data points from sensors instead of simulated data would \textcolor{black}{result in higher heterogeneity in the data and ensure robustness in the models}. The constrained eco-routing concept is an important aspect to be tackled and illustrates the trade offs compared to the regular eco-routing application. As 100\% CAVs MPR is employed in this study, the impact of different MPRs should be taken into account. The most preferable MPR of CAVs could vary from a traffic condition to another and this has to be defined. The predictive models utilized in this study can be further enhanced by using more data points to represent the conditions of low frequency of occurrence. In addition, predictive models can be developed based on categorized characteristics e.g. speed limit, number of lanes, etc. of links for further enhancements. With regards to the scaleability aspect, it is suggested that anticipatory routing is applied in a larger network with both uninterrupted and interrupted traffic flow. Nevertheless, the predictive models should accommodate the difference between the two types of traffic flow. The employed distributed routing framework adopts only the V2I and I2I communication. The impact of incorporating the V2V communication is another suggestion for the future work. Incorporating anticipatory routing strategies as options in the personal navigation platforms would contribute to more efficient and sustainable transportation systems. Despite the strong contradiction between the NOx and speed, including the NOx as a routing objective is preferred.}





\bibliographystyle{abbrv}
\bibliography{references.bib}
\end{document}